\documentclass[11pt]{article}%
\usepackage{amsmath}
\usepackage{amsfonts}
\usepackage{amssymb}
\usepackage{graphicx}
\usepackage{theorem}
\usepackage{makeidx}
\usepackage{verbatim}
\usepackage{color}%
\providecommand{\U}[1]{\protect\rule{.1in}{.1in}}
\newtheorem{thm}{Theorem}[section]
\newtheorem{lm}[thm]{Lemma}
\newtheorem{pr}[thm]{Proposition}
\newtheorem{df}[thm]{Definition}
\newtheorem{rmk}[thm]{Remark}
\newtheorem{cor}[thm]{Corollary}
{\theorembodyfont{\upshape}
\newtheorem{examp}[thm]{Example}
}
\numberwithin{equation}{section} 
\begin{document}

\title{ }

\begin{center}
\vspace*{1.3cm}

\textbf{ON SUBREGULARITY PROPERTIES\ OF\ SET-VALUED\ MAPPINGS.\ APPLICATIONS
TO SOLID\ VECTOR\ OPTIMIZATION}

\bigskip

by

\bigskip

MARIUS APETRII\footnote{{\small {Faculty of Mathematics, "Al. I. Cuza"
University,} \ {Bd. Carol I, nr. 11, 700506 -- Ia\c{s}i, Romania, e-mail:
\texttt{mapetrii@uaic.ro}}}}, MARIUS DUREA\footnote{{\small {Faculty of
Mathematics, "Al. I. Cuza" University,} \ {Bd. Carol I, nr. 11, 700506 --
Ia\c{s}i, Romania, e-mail: \texttt{durea@uaic.ro}}}} and
RADU\ STRUGARIU\footnote{{\small Department of Mathematics, "Gh. Asachi"
Technical University, \ {Bd. Carol I, nr. 11, 700506 -- Ia\c{s}i, Romania,
e-mail: \texttt{rstrugariu@tuiasi.ro}}}}

\end{center}

\bigskip

\bigskip

\noindent{\small {\textbf{Abstract:}} In this work we classify the at-point
regularities of set-valued mappings into two categories and then we analyze
their relationship through several implications and examples. After this
theoretical tour, we use the subregularity properties to deduce implicit
theorems for set-valued maps. Finally, we present some applications to the
study of multicriteria optimization problems.}

\bigskip


\bigskip


\section{Introduction}

This work is a natural continuation of the recent papers \cite{DurStr4} and
\cite{DNS2011}. We mention that in \cite{DurStr4} several results concerning
the behavior of the implicit solution mappings associated to parametric
variational systems are given under around-point regularity assumptions for
the initial mappings, while in \cite{DNS2011} the same concepts are employed
in order to get necessary optimality conditions for multicriteria
optimization. However, let us remark that, on one hand, the calmness property
replaces the Aubin continuity in many recent works in literature in order to
get well-posedness results for parametric systems (see, for instance,
\cite{ArtMord2010}, \cite{CKY}, \cite{LiMor} and the references therein) and,
on the other hand, in certain situations, metric subregularity of the
constraint system is enough for deriving necessary optimality conditions for
generalized mathematical programs (see \cite{YY}, \cite{HO2005}).

In this perspective, we revisit some results in \cite{DurStr4}, \cite{DNS2011}%
, trying to replace, where possible, the around-point regularities of
set-valued mappings by weaker concepts of at-point regularity. First of all,
we classify the regularity notions at the reference point into two categories,
which we call first and, respectively, second type at-point regularity.
Following the pattern stated for regularity around the reference point, for
both these types we consider the corresponding triads of openness, metric
regularity and Lipschitzness. For the first type we have the at-point openness
at linear rate (\cite{Urs1996}), the metric hemiregularity (\cite{ArtMord2010}%
) and the pseudocalmness (\cite{DurStr4}). For the second type, which seems to
be of greater interest compared to the first one, up to now there exist only
two concepts: the calmness and the metric subregularity (see \cite[Section
3H]{DontRock2009b}). Our first aim is to complete the triad by an appropriate
equivalent openness notion (which we call the linear pseudo-openness), giving
us the possibility to have a deeper insight on the results involving the
second type of regularity at the reference point. After that, we establish the
relationship between all these notions, with a special emphasis on the case of
linear bounded operators.

With all these facts in mind, we are able to discuss our main results, which
are divided into three themes. The first one concerns implicit multifunction
type theorems using at-point regularity of the second type, showing the
technical advantages of the linear pseudo-openness concept. This is in the
line with the remark that, in general, the linear openness notions are
technically easier to deal with in the proofs, while the equivalent metric
regularities are more useful in applications (see \cite{DMO}, \cite{Dmi2005}).
The second theme is devoted to the analysis of the local-sum stability, a
notion recently introduced in \cite{DurStr4} in relation with the conservation
of the Aubin property for sum-type set-valued maps. In this work we study the
natural link between this sort of stability and the conservation of calmness
at summation, and we emphasize by examples other situations where the sum
multifunction is calm, in the absence of this property for the component
mappings. The third theme is dedicated to the study of variational systems in
the context of the second type at-point regularity, taking advantage of the
technical analysis developed before. Firstly, we derive a result concerning
the metric subregularity of the implicit mappings associated to parametric
variational systems, and we provide an example which shows that the sole
around-point regularity assumption used to this aim cannot be dropped.
Secondly, we provide a fixed-point result for the parametric case of
composition of two set-valued mappings, which can be seen in relation to some
recent metric extensions of the Lyusternik-Graves theorem (see
\cite{DonFra2010}, \cite{DonFra2011}, \cite{DurStr5}, \cite{DurStr6}). Then we
use the above mentioned result to deduce, under additional appropriate
conditions, the calmness of the solution mappings of the parametric systems.

In the final section of the paper we employ the second type of at-point
regularity in the study of solid set-valued optimization problems. First of
all, using the metric subregularity of the constraint system, we combine a
Clarke type penalization technique and a scalarization result in order to
reduce the problem of getting necessary optimality conditions for weak Pareto
minimizers to the one of finding local minimizers for a scalar function. After
that, we use some ideas from \cite{DNS2011} in order to deduce sufficient
conditions for the needed metric subregularity of the constraint system,
following the error bounds approach and using the Mordukhovich generalized
differentiation objects for the formulation of regularity by means of dual
objects. Putting all the previous facts together, we finally get the expected
necessary optimality conditions in terms of Mordukhovich differentiation, by
expressing some generalized Lagrange multipliers rules in the normal form for
the proposed optimization problem.

\section{Preliminaries}

In what follows, we suppose that the involved spaces are metric spaces, unless
otherwise stated. In this setting, $B(x,r)$ and $D(x,r)$ denote the open and
the closed ball with center $x$ and radius $r,$ respectively. On a product
space we take the additive metric. If $x\in X$ and $A\subset X,$ one defines
the distance from $x$ to $A$ as $d(x,A):=\inf\{d(x,a)\mid a\in A\}.$ As usual,
we use the convention $d(x,\emptyset)=\infty.$ The excess from a set $A$ to a
set $B$ is defined as $e(A,B):=\sup\{d(a,B)\mid a\in A\}.$ For a non-empty set
$A\subset X$ we put $\operatorname*{cl}A$ for its topological closure. One
says that a set $A$ is locally closed around $x\in A$ if there exists $r>0$
such that $A\cap D(x,r)$ is closed.

Let $F:X\rightrightarrows Y$ be a multifunction. The domain and the graph of
$F$ are denoted respectively by $\operatorname*{Dom}F:=\{x\in X\mid
F(x)\neq\emptyset\}$ and $\operatorname*{Gr}F:=\{(x,y)\in X\times Y\mid y\in
F(x)\}.$ If $A\subset X$ then $F(A):=%
{\displaystyle\bigcup\limits_{x\in A}}
F(x).$ The inverse set-valued map of $F$ is $F^{-1}:Y\rightrightarrows X$
given by $F^{-1}(y)=\{x\in X\mid y\in F(x)\}$.

We denote by $P$ the metric space of parameters. For a (parametric)
multifunction $F:X\times P\rightrightarrows Y,$ we use the notations:
$F_{p}(\cdot):=F(\cdot,p)$ and $F_{x}(\cdot):=F(x,\cdot).$

We divide the reminder of this section into three different subsections, each
one being dedicated to a certain type of regularity in set-valued setting.

\subsection{Around-point regularity}

We recall now the concepts of openness at linear rate, metric regularity and
Aubin property of a multifunction around the reference point. Generally, when
one speaks about regularity for a set-valued map, one refers to these concepts.

\begin{df}
\label{around}Let $L>0,$ $F:X\rightrightarrows Y$ be a multifunction and
$(\overline{x},\overline{y})\in\operatorname{Gr}F.$

(i) $F$ is said to be open at linear rate $L,$ or $L-$open around
$(\overline{x},\overline{y}),$ if there exist a positive number $\varepsilon
>0$ and two neighborhoods $U\in\mathcal{V}(\overline{x}),$ $V\in
\mathcal{V}(\overline{y})$ such that, for every $\rho\in(0,\varepsilon)$ and
every $(x,y)\in\operatorname*{Gr}F\cap\lbrack U\times V],$%
\begin{equation}
B(y,\rho L)\subset F(B(x,\rho)). \label{Lopen}%
\end{equation}

The supremum of $L>0$ over all the combinations $(L,U,V,\varepsilon)$ for
which (\ref{Lopen}) holds is denoted by $\operatorname*{lop}F(\overline
{x},\overline{y})$ and is called the exact linear openness bound, or the exact
covering bound of $F$ around $(\overline{x},\overline{y}).$

(ii) $F$ is said to have the Aubin property around $(\overline{x},\overline
{y})$ with constant $L$ if there exist two neighborhoods $U\in\mathcal{V}%
(\overline{x}),$ $V\in\mathcal{V}(\overline{y})$ such that, for every $x,u\in
U,$%
\begin{equation}
e(F(x)\cap V,F(u))\leq Ld(x,u). \label{LLip_like}%
\end{equation}

The infimum of $L>0$ over all the combinations $(L,U,V)$ for which
(\ref{LLip_like}) holds is denoted by $\operatorname*{lip}F(\overline
{x},\overline{y})$ and is called the exact Lipschitz bound of $F$ around
$(\overline{x},\overline{y}).$

(iii) $F$ is said to be metrically regular around $(\overline{x},\overline
{y})$ with constant $L$ if there exist two neighborhoods $U\in\mathcal{V}%
(\overline{x}),$ $V\in\mathcal{V}(\overline{y})$ such that, for every
$(x,y)\in U\times V,$%
\begin{equation}
d(x,F^{-1}(y))\leq Ld(y,F(x)). \label{Lmreg}%
\end{equation}

The infimum of $L>0$ over all the combinations $(L,U,V)$ for which
(\ref{Lmreg}) holds is denoted by $\operatorname*{reg}F(\overline{x}%
,\overline{y})$ and is called the exact regularity bound of $F$ around
$(\overline{x},\overline{y}).$
\end{df}

The next proposition contains the well-known links between the notions
presented above. See \cite{RocWet}, \cite{Mor2006}, \cite{DontRock2009b} for
more details about its proof and for historical facts.

\begin{pr}
\label{link_around}Let $F:X\rightrightarrows Y$ be a multifunction and
$(\overline{x},\overline{y})\in\operatorname{Gr}F.$ Then $F$ is open at linear
rate around $(\overline{x},\overline{y})$ iff $F^{-1}$ has the Aubin property
around $(\overline{y},\overline{x})$ iff $F$ is metrically regular around
$(\overline{x},\overline{y})$. Moreover, in every of the previous situations,%
\[
(\operatorname*{lop}F(\overline{x},\overline{y}))^{-1}=\operatorname*{lip}%
F^{-1}(\overline{y},\overline{x})=\operatorname*{reg}F(\overline{x}%
,\overline{y}).
\]

\end{pr}

All the regularity concepts given before have parametric counterparts, which
we present next.

\begin{df}
Let $L>0,$ $F:X\times P\rightrightarrows Y$ be a multifunction and
$((\overline{x},\overline{p}),\overline{y})\in\operatorname{Gr}F.$

(i) $F$ is said to be open at linear rate $L,$ or $L-$open, with respect to
$x$ uniformly in $p$ around $((\overline{x},\overline{p}),\overline{y})$ if
there exist a positive number $\varepsilon>0$ and some neighborhoods
$U\in\mathcal{V}(\overline{x}),$ $V\in\mathcal{V}(\overline{p}),$
$W\in\mathcal{V}(\overline{y})$ such that, for every $\rho\in(0,\varepsilon),$
every $p\in V$ and every $(x,y)\in\operatorname*{Gr}F_{p}\cap\lbrack U\times
W],$%
\begin{equation}
B(y,\rho L)\subset F_{p}(B(x,\rho)). \label{pLopen}%
\end{equation}

The supremum of $L>0$ over all the combinations $(L,U,V,W,\varepsilon)$ for
which (\ref{pLopen}) holds is denoted by $\widehat{\operatorname*{lop}}%
_{x}F((\overline{x},\overline{p}),\overline{y})$ and is called the exact
linear openness bound, or the exact covering bound of $F$ in $x$ around
$((\overline{x},\overline{p}),\overline{y}).$

(ii) $F$ is said to have the Aubin property with respect to $x$ uniformly in
$p$ around $((\overline{x},\overline{p}),\overline{y})$ with constant $L$ if
there exist some neighborhoods $U\in\mathcal{V}(\overline{x}),$ $V\in
\mathcal{V}(\overline{p}),$ $W\in\mathcal{V}(\overline{y})$ such that, for
every $x,u\in U$ and every $p\in V,$%
\begin{equation}
e(F_{p}(x)\cap W,F_{p}(u))\leq Ld(x,u). \label{pLLip_like}%
\end{equation}

The infimum of $L>0$ over all the combinations $(L,U,V,W)$ for which
(\ref{pLLip_like}) holds is denoted by $\widehat{\operatorname*{lip}}%
_{x}F((\overline{x},\overline{p}),\overline{y})$ and is called the exact
Lipschitz bound of $F$ in $x$ around $((\overline{x},\overline{p}%
),\overline{y}).$

(iii) $F$ is said to be metrically regular with respect to $x$ uniformly in
$p$ around $((\overline{x},\overline{p}),\overline{y})$ with constant $L$ if
there exist some neighborhoods $U\in\mathcal{V}(\overline{x}),$ $V\in
\mathcal{V}(\overline{p}),$ $W\in\mathcal{V}(\overline{y})$ such that, for
every $(x,p,y)\in U\times V\times W,$%
\begin{equation}
d(x,F_{p}^{-1}(y))\leq Ld(y,F_{p}(x)). \label{pLmreg}%
\end{equation}

The infimum of $L>0$ over all the combinations $(L,U,V,W)$ for which
(\ref{pLmreg}) holds is denoted by $\widehat{\operatorname*{reg}}%
_{x}F((\overline{x},\overline{p}),\overline{y})$ and is called the exact
regularity bound of $F$ in $x$ around $((\overline{x},\overline{p}%
),\overline{y}).$
\end{df}

The corresponding notions with respect to $p$ uniformly in $x$ can be written similarly.

In the sequel, we emphasize the fact that the corresponding "at-point"
properties could be separated into two different categories, which for the
sake of clarity we present as type I and type II.

\subsection{At-point regularity: type I}

The first type of at-point regularity contains the linear openness at point,
the pseudocalmness and the metric hemiregularity, as follows.

\begin{df}
\label{at}Let $L>0,$ $F:X\rightrightarrows Y$ be a multifunction and
$(\overline{x},\overline{y})\in\operatorname{Gr}F.$

(i) $F$ is said to be open at linear rate $L,$ or $L-$open at $(\overline
{x},\overline{y})$ if there exists a positive number $\varepsilon>0$ such
that, for every $\rho\in(0,\varepsilon),$%
\begin{equation}
B(\overline{y},\rho L)\subset F(B(\overline{x},\rho)). \label{Lopen_at}%
\end{equation}

The supremum of $L>0$ over all the combinations $(L,\varepsilon)$ for which
(\ref{Lopen_at}) holds is denoted by $\operatorname*{plop}F(\overline
{x},\overline{y})$ and is called the exact punctual linear openness bound of
$F$ at $(\overline{x},\overline{y}).$

(ii) $F$ is said to be pseudocalm with constant $L,$ or $L-$pseudocalm at
$(\overline{x},\overline{y}),$ if there exists a neighborhood $U\in
\mathcal{V}(\overline{x})$ such that, for every $x\in U,$%
\begin{equation}
d(\overline{y},F(x))\leq Ld(x,\overline{x}). \label{Lpseudocalm}%
\end{equation}

The infimum of $L>0$ over all the combinations $(L,U)$ for which
(\ref{Lpseudocalm}) holds is denoted by $\operatorname*{psdclm}F(\overline
{x},\overline{y})$ and is called the exact bound of pseudocalmness for $F$ at
$(\overline{x},\overline{y}).$

(iii) $F$ is said to be metrically hemiregular with constant $L,$ or
$L-$metrically hemiregular at $(\overline{x},\overline{y}),$ if there exists a
neighborhood $V\in\mathcal{V}({\overline{y}})$ such that, for every $y\in V,$%
\begin{equation}
d(\overline{x},F^{-1}(y))\leq Ld(y,{\overline{y})}. \label{Lhemireg}%
\end{equation}

The infimum of $L>0$ over all the combinations $(L,V)$ for which
(\ref{Lhemireg}) holds is denoted by $\operatorname*{hemreg}F(\overline
{x},\overline{y})$ and is called the exact hemiregularity bound of $F$ at
$(\overline{x},\overline{y}).$
\end{df}

As in the case of around-point regularity, some equivalences between these
notions hold. For the proof, see, for instance, \cite[Proposition
2.4]{DurStr4}.

\begin{pr}
\label{link_at}Let $L>0,$ $F:X\rightrightarrows Y$ and $(\overline
{x},\overline{y})\in\operatorname{Gr}F.$ Then $F$ is $L-$open at
$(\overline{x},\overline{y})$ iff $F^{-1}$ is $L^{-1}-$pseudocalm at
$(\overline{y},\overline{x})$ iff $F$ is $L^{-1}-$metrically hemiregular at
$(\overline{x},\overline{y})$. Moreover, in every of the previous situations,%
\[
(\operatorname*{plop}F(\overline{x},{\overline{y}}))^{-1}%
=\operatorname*{psdclm}F^{-1}({\overline{y},}\overline{x}%
)=\operatorname*{hemreg}F(\overline{x},{\overline{y}}).
\]

\end{pr}

Let us mention that partial corresponding variants could be done as in the
previous subsection.

\subsection{At-point regularity: type II}

The second type of at-point regularity contains the calmness, the metric
subregularity, and a notion we introduce here under the name of linear
pseudo-openness. This novelty serves to complete the regularity triad in this
case. Moreover, it proves to be useful in the attempt of getting implicit
multifunction results under weaker assumptions, and this is the main aim of
this work.

\begin{df}
\label{at_2}Let $L>0,$ $F:X\rightrightarrows Y$ be a multifunction and
$(\overline{x},\overline{y})\in\operatorname{Gr}F.$

(i) $F$ is said to be linearly pseudo-open with modulus $L,$ or $L-$%
pseudo-open at $(\overline{x},\overline{y})$ if there exist $U\in
\mathcal{V}(\overline{x})$ and $\varepsilon>0$ such that for every $\rho
\in(0,\varepsilon)$ and for every $x\in U\cap F^{-1}(B(\overline{y},L\rho)),$%
\begin{equation}
\overline{y}\in F(B(x,\rho)). \label{psdop}%
\end{equation}

The supremum of $L>0$ over all the combinations $(L,U,\varepsilon)$ for which
(\ref{psdop}) holds is denoted by $\operatorname*{lpo}F(\overline{x}%
,\overline{y})$ and is called the exact linear pseudo-openness bound of $F$ at
$(\overline{x},\overline{y}).$

(ii) $F$ is said to be calm with constant $L,$ or $L-$calm at $(\overline
{x},\overline{y}),$ if there exists some neighborhoods $U\in\mathcal{V}%
(\overline{x}),$ $V\in\mathcal{V}(\overline{y})$ such that, for every $x\in
U,$%
\begin{equation}
e(F(x)\cap V,F(\overline{x}))\leq Ld(x,\overline{x}). \label{calm}%
\end{equation}

The infimum of $L>0$ over all the combinations $(L,U,V)$ for which
(\ref{calm}) holds is denoted by $\operatorname*{clm}F(\overline{x}%
,\overline{y})$ and is called the exact bound of calmness for $F$ at
$(\overline{x},\overline{y}).$

(iii) $F$ is said to be metrically subregular with constant $L,$ or
$L-$metrically subregular at $(\overline{x},\overline{y})$ if there exists a
neighborhood $U\in\mathcal{V}(\overline{x})$ such that, for every $x\in U,$%
\begin{equation}
d(x,F^{-1}(\overline{y}))\leq Ld(\overline{y},F(x)). \label{subreg}%
\end{equation}

The infimum of $L>0$ over all the combinations $(L,U)$ for which
(\ref{Lhemireg}) holds is denoted by $\operatorname*{subreg}F(\overline
{x},\overline{y})$ and is called the exact subregularity bound of $F$ at
$(\overline{x},\overline{y}).$
\end{df}

\bigskip

Remark that the properties of calmness and metrically subregularity are well
known. The first concept of Definition \ref{at_2} which we introduce here
plays the same role in this triplet as the openness in the above well-known
ones. The exact meaning of this assertion is given in the next result.

\begin{pr}
\label{link_at2}Let $F:X\rightrightarrows Y$ and $(\overline{x},\overline
{y})\in\operatorname{Gr}F.$ Then $F$ is linearly pseudo-open at $(\overline
{x},\overline{y})$ iff $F^{-1}$ is calm at $(\overline{y},\overline{x})$ iff
$F$ is metrically subregular at $(\overline{x},\overline{y})$. Moreover, in
every of the previous situations,%
\[
(\operatorname*{lpo}F(\overline{x},{\overline{y}}))^{-1}=\operatorname*{clm}%
F^{-1}({\overline{y},}\overline{x})=\operatorname*{subreg}F(\overline
{x},{\overline{y}}).
\]

\end{pr}

\noindent\textbf{Proof.} The equivalence between the calmness of $F^{-1}$ and
the metric subregularity of $F$, as well as the relation between the
corresponding regularity moduli are well-known (see, for instance,
\cite[Section 3H]{DontRock2009b}).

Let us prove the equivalence between the linear pseudo-openness of $F$ and the
calmness of $F^{-1}.$ Suppose first that $F$ is linearly pseudo-open at
$(\overline{x},\overline{y})$ with modulus $L>0.$ Then there exist
$U\in\mathcal{V}(\overline{x})$ and $\varepsilon>0$ such that, for every
$\rho\in(0,\varepsilon)$ and every $x\in U\cap F^{-1}(B(\overline{y},L\rho)),$
$\overline{y}\in F(B(x,\rho)).$ Consider $\varepsilon^{\prime}:=L\varepsilon
,V:=B(\overline{y},\varepsilon^{\prime}),$ and take $y\in V$ and $x\in
F^{-1}(y)\cap U.$ Without loosing the generality, suppose that $y\not =%
\overline{y},$ because otherwise the desired relation trivially holds. Then,
there exists $\rho\in(0,\varepsilon)$ such that $d(y,\overline{y}%
)=L\rho<L\varepsilon.$ Take $\tau>0$ arbitrary small such that $\rho^{\prime
}:=\rho+\tau<\varepsilon$. Consequently, $x\in U\cap F^{-1}(B(\overline
{y},L\rho^{\prime})),$ hence $\overline{y}\in F(B(x,\rho^{\prime})),$ or there
exists $z\in B(x,\rho^{\prime})$ such that $z\in F(\overline{y}).$ But this
means that
\[
d(x,F(\overline{y}))\leq d(x,z)\leq\rho+\tau=L^{-1}d(y,\overline{y})+\tau.
\]
Because $x$ was arbitrarily taken from $F^{-1}(y)\cap U,$ from the previous
relation one deduces that
\[
e(F^{-1}(y)\cap U,F(\overline{y}))\leq L^{-1}d(y,\overline{y})+\tau
\]

\noindent for every $\tau>0$ arbitrary small, so making $\tau\rightarrow0$ one
deduces that $F^{-1}$ is calm at $(\overline{y},\overline{x})$ with modulus
$L^{-1}.$

Suppose now that $F^{-1}$ is calm at $(\overline{y},\overline{x})$ with
modulus $L^{-1},$ so there exists $U\in\mathcal{V}(\overline{x})$ and
$V\in\mathcal{V}(\overline{y})$ such that, for every $y\in V,$%
\[
e(F^{-1}(y)\cap U,F(\overline{y}))\leq L^{-1}d(y,\overline{y}).
\]

We will prove that $F$ is linearly pseudo-open at $(\overline{x},\overline
{y})$ with modulus smaller, but arbitrarily close to $L.$ Choose $\xi>0$ such
that $L-\xi>0,$ take $\varepsilon>0$ such that $B(\overline{y},L\varepsilon
)\subset V,$ fix arbitrary $\rho\in(0,\varepsilon)$ and $x\in U\cap
F^{-1}(B(\overline{y},(L-\xi)\rho)).$ Then there exists $y\in B(\overline
{y},(L-\xi)\rho)\subset B(\overline{y},L\rho)\subset B(\overline
{y},L\varepsilon)\subset V$ such that $x\in F^{-1}(y).$ Therefore,
$d(x,F^{-1}(\overline{y}))\leq L^{-1}d(y,\overline{y}),$ so one can find $z\in
F^{-1}(\overline{y})$ such that%
\[
d(x,z)<L^{-1}d(y,\overline{y})+L^{-1}\rho\xi<L^{-1}(L-\xi)\rho+L^{-1}\rho
\xi=\rho.
\]

\noindent In conclusion, there exists $z\in B(x,\rho)$ such that $\overline
{y}\in F(z),$ i.e. the conclusion.$\hfill\square$

\bigskip

\begin{rmk}
\label{char}Notice that, in our notation, we obviously have the following characterizations:

(i) $F$ is calm at\ $(\overline{x},\overline{y})$ with constant $L>0$ if and
only if there exist $U\in\mathcal{V}(\overline{x})$ and $V\in\mathcal{V}%
(\overline{y})$ such that for every $x\in U$ and $y\in V$ with
$d(y,F(\overline{x}))>Ld(x,\overline{x}),\ $we have $y\notin F(x);$

(ii) $F$ is metrically subregular at\ $(\overline{x},\overline{y})$ with
constant $L>0$ if and only if there exist $U\in\mathcal{V}(\overline{x})$ and
$V\in\mathcal{V}(\overline{y})$ such that for every $x\in U$ and $y\in V$ with
$d(x,F^{-1}(\overline{y}))>Ld(y,\overline{y}),$ we have $y\notin F(x).$
\end{rmk}

\bigskip

Even if one can define all the corresponding partial notions to the concepts
in Definition \ref{at_2}, we restrict ourselves to the case of calmness,
because this is exactly what we use in the next sections. More precisely, $F$
is said to be calm with respect to $x$ uniformly in $p$ at $((\overline
{x},\overline{p}),\overline{y})$ with constant $L>0$ if there exist some
neighborhoods $U\in\mathcal{V}(\overline{x}),$ $V\in\mathcal{V}(\overline
{p}),$ $W\in\mathcal{V}(\overline{y})$ such that, for every $x\in U$ and every
$p\in V,$%
\begin{equation}
e(F_{p}(x)\cap W,F_{p}(\overline{x}))\leq Ld(x,\overline{x}). \label{clmFxp}%
\end{equation}

The infimum of $L>0$ over all the combinations $(L,U,V,W)$ for which
(\ref{clmFxp}) holds is denoted by $\widehat{\operatorname*{clm}}%
_{x}F((\overline{x},\overline{p}),\overline{y})$ and is called the exact
calmness bound of $F$ in $x$ at $((\overline{x},\overline{p}),\overline{y}).$

\bigskip

We close this section by some comments concerning the comparison of the three
sets of concepts we have listed above. It is well known (and easy to see) that
around-point regularity implies both types of at-point regularity while the
converse implications obviously fail.

A remarkable situations where several interesting assertions could be
additionally established is that of linear continuous operators, i.e. the case
where $F$ is replaced by $A\in\mathcal{L}(X,Y),$ where $\mathcal{L}(X,Y)$
denotes the normed vector space of linear bounded operators acting between $X$
and $Y$. A consequence of Banach Open Principle is that $A$ is open at linear
rate around a (every) $(x,Ax)$ if and only if $A$ is surjective. Moreover,
following \cite[Proposition 5.2]{ArtMord2010}, this is further equivalent to
the hemiregularity of $A$ at a (every) $(x,Ax)$. Since for a linear continuous
operator all the regularity moduli are not depending on the reference point we
can remove it from their notations and, under surjectivity of $A,$ one has
\begin{equation}
\operatorname*{psdclm}A^{-1}=\operatorname*{lip}A^{-1}=\operatorname*{hemreg}%
A=\operatorname*{reg}A=(\operatorname*{plop}A)^{-1}=(\operatorname*{lop}%
A)^{-1}=\left\Vert (A^{\ast})^{-1}\right\Vert , \label{mod}%
\end{equation}
\noindent where $A^{\ast}\in\mathcal{L}(Y^{\ast},X^{\ast})$ denotes the
adjoint operator of $A.$

For the case of second type at-point regularities the situation changes
significantly. First of all, if $X$ and $Y$ are finite dimensional, then every
$A\in\mathcal{L}(X,Y)$ is metrically subregular at every point of its graph.
To see this, observe that if $X,Y$ are finite dimensional and $A\in
\mathcal{L}(X,Y),$ then $\operatorname*{Im}A$ is isomorphic with $X_{2},$
where by $X_{2}$ we denote the algebraic complement of $\operatorname*{Ker}A$
in $X.$ One has to prove that there is $L>0$ such that, for every $x\in X,$%
\begin{equation}
\inf\{\left\Vert x-u\right\Vert \mid u\in\operatorname*{Ker}A\}\leq
L\left\Vert Ax\right\Vert . \label{ineg_fd}%
\end{equation}

Observe first that if $x\in\operatorname*{Ker}A,$ then the above inequality
trivially holds. Define next the isomorphism $A_{1}:X_{2}\rightarrow
\operatorname*{Im}A,$ given by $A_{1}x:=Ax$ for every $x\in X_{2}$, and apply
the Banach Open Principle for $A_{1}$ to deduce that $A_{1}$ is open.
Equivalently, there exists $L>0$ such that, for every $y\in\operatorname*{Im}%
A,$ there exists $x\in X_{2}$ such that $y=A_{1}x$ and $\left\Vert
x\right\Vert \leq L\left\Vert y\right\Vert .$ Take arbitrary $x\in X_{2}.$
Then it uniquely corresponds to $y=A_{1}x=Ax\in\operatorname*{Im}A,$ hence%
\[
\left\Vert x\right\Vert \leq L\left\Vert Ax\right\Vert .
\]

Finally, for an arbitrary $x\in X,$ we decompose it as $x=x^{\prime}%
+x^{\prime\prime},$ with $x^{\prime}\in\operatorname*{Ker}A$ and
$x^{\prime\prime}\in X_{2},$ and we have%
\[
\inf\{\left\Vert x-u\right\Vert \mid u\in\operatorname*{Ker}A\}\leq\left\Vert
x-x^{\prime}\right\Vert =\left\Vert x^{\prime\prime}\right\Vert \leq
L\left\Vert Ax^{\prime\prime}\right\Vert =L\left\Vert Ax\right\Vert ,
\]

\noindent which ends the proof. Note that an alternative proof can be given as
follows: one knows that the distance in the left-hand side of (\ref{ineg_fd})
is attained and can be written as $\left\langle x,x^{\ast}\right\rangle $
where $x^{\ast}$ belongs to orthogonal subspace of $\operatorname*{Ker}A$ and
its operatorial norm is $1$ (see \cite[Theorem 3.8.4 (vii)]{Zal2002}). By
means of Farkas Lemma (\cite[Corrolary 22.3.1]{Roc1970}), $x^{\ast}$ can be
expressed as a linear combination of the scalar linear components mappings of
$A$ and then the thesis is proved by a manipulation of some easy inequalities.
For other details concerning these aspects, see the comments in \cite[Section
3H]{DontRock2009b}. The previous discussion means that any linear operator on
finite dimensional spaces which is not surjective is metrically subregular,
but fails to be metrically regular. Therefore, even for linear operators,
metric subregularity does not imply metric regularity.

On infinite dimensional spaces, there exist linear bounded operators which
fail to be metrically subregular. For instance, consider the spaces $m$ and
$l^{2}$ of bounded and respectively square-summable sequences of real numbers
with their usual norms and $T:m\rightarrow l^{2}$ with $T((x_{n}%
))=(n^{-1}x_{n})$ for every sequence $(x_{n})_{n\in\mathbb{N}\setminus\{0\}}.$
It is easy to show that $T$ is well defined, linear, continuous $(\left\Vert
T\right\Vert =\sqrt{6^{-1}}\pi)$ and injective. Nevertheless, supposing that
$T$ would be metrically subregular at $(0,0),$ then should exist $L>0$ s.t.
for every $(x_{n})\in m,$%
\[
\left\Vert (x_{n})\right\Vert \leq L\left\Vert T(x_{n})\right\Vert .
\]
Taking, for every natural $k\neq0,$ $(x_{n}^{k})_{n}\ $as the sequence having
all the components zero except that on $k-$th place which is $1,$ then the
above relation reads as $1\leq Lk^{-1}$ for every $k\in\mathbb{N}%
\setminus\{0\}$. This is a contradiction which can be eliminated only if $T$
is not metrically subregular at $(0,0).$

On the other hand, on Banach spaces, if $A$ is injective and its image is
closed, then this is equivalent to the following property:
\begin{equation}
\exists M>0\text{ s.t. }\forall x\in X,\text{ }\left\Vert x\right\Vert \leq
M\left\Vert Ax\right\Vert . \label{mg_inf}%
\end{equation}
It is easy to observe that this last relation easily leads to the linear
pseudo-openness at every point of the graph. Indeed, it is enough to see that
in these assumptions one can write for $(\overline{x},\overline{y})=(0,0)$ and
every $x\in X$%
\[
d(x,A^{-1}(0))\leq\left\Vert x-0\right\Vert \leq M\left\Vert Ax\right\Vert
=Md(0,Ax),
\]
which ensures the metric subregularity of $A.$ With this remark, we infer that
in the above example, the metric subregularity fails because the image of $T$
is not closed in $l^{2}.$

\bigskip

Let us observe now that if $A$ is open (hence surjective) with constant $L>0$
then obviously $A$ is pseudo-open with modulus $L$ and $\operatorname*{lop}%
A\leq\operatorname*{lpo}A.$

The opposite inequality is also true. Indeed, if $A$ is pseudo-open then we
find $\varepsilon>0$ and a neighborhood $U\in\mathcal{V}(0)$ such that for
every $\rho\in(0,\varepsilon)$ and for every $x\in U\cap A^{-1}(B(0,L\rho))$
we have that $0\in A(B(x,\rho)).$ Starting with $y\in B(0,L\rho),$ since $A$
is surjective, then there exists $x\in X$ with $Ax=y,$ and taking into account
that $U$ is a neighborhood of $0$ then we find $\lambda\geq1$ and $x^{\prime
}\in U$ such that $x=\lambda x^{\prime}.$ Now, we have that $\left\Vert
Ax^{\prime}\right\Vert \leq\lambda^{-1}\rho L,$ therefore\ $0\in
A(B(x^{\prime},\lambda^{-1}\rho))$ because $x^{\prime}\in U\cap A^{-1}%
(B(0,\lambda^{-1}\rho L)).$ From this and from linearity of $A$ we obtain that
$0\in A(B(x,\rho)),$ which shows that $A$ is open at linear rate $L$ and,
taking into account (\ref{mod}), $\operatorname*{lpo}A\leq\operatorname*{plop}%
A=\operatorname*{lop}A.$ In conclusion, under surjectivity of $A,$ one can add
to (\ref{mod}) the following chain of equalities:%
\[
\operatorname*{clm}A^{-1}=\operatorname*{lip}A^{-1}=\operatorname*{subreg}%
A=\operatorname*{reg}A=(\operatorname*{lpo}A)^{-1}=(\operatorname*{lop}%
A)^{-1}=\left\Vert (A^{\ast})^{-1}\right\Vert .
\]

\section{Main results}

This section is divided into three subsections, as follows:

\begin{itemize}
\item the first one consists of an implicit multifunction type theorem
displaying at-point regularity of the second type;

\item the second one gives further insights on recently introduced notion of
local sum-stability of two multifunctions; more precisely, we investigate the
relation with the calmness of the sum of two mappings both through theoretical
results and examples;

\item the third one concerns the study of parametric variational systems in
the context of the second type at-point regularity.
\end{itemize}

\subsection{At-point regularity for implicit multifunctions}

In this subsection we obtain a result concerning at-point regularity of the
second type for a generalized implicit set-valued map.

To this aim, we consider a setting coming from the study of parametric
variational systems. Given a multifunction $H:X\times P\rightrightarrows Y,$
define the implicit mapping $S:P\rightrightarrows X$ by%
\begin{equation}
S(p):=\{x\in X\mid0\in H(x,p)\}. \label{S}%
\end{equation}

\noindent The study of the well-posedness properties of $S$, including the
so-called Robinson regularity, was a constant issue of operations research in
the last four decades (see \cite{DontRock2009b} for a extended discusion and
historical comments).

The next result is in the line of \cite[Theorem 3.6]{DurStr4}, but for the
newly introduced concept of linear pseudo-openness instead of the classical
genuine linear openness around the reference point, and for calmness instead
of Aubin property.

\begin{thm}
\label{main}Let $X,P$ be metric spaces, $Y$ be a normed vector space,
$H:X\times P\rightrightarrows Y$ be a set-valued map and $(\overline
{x},\overline{p},0)\in\operatorname*{Gr}H.$ Denote by $H_{p}(\cdot
):=H(\cdot,p),$ $H_{x}(\cdot):=H(x,\cdot).$

(i) If $H_{\overline{p}}$ is linearly pseudo-open with modulus $c>0$ at
$(\overline{x},0),$ then there exist $\alpha,\beta>0$ such that, for every
$x\in B(\overline{x},\alpha),$%
\begin{equation}
d(x,S(\overline{p}))\leq c^{-1}d(0,H(x,\overline{p})\cap B(0,\beta)).
\label{clmS}%
\end{equation}

If, moreover, $H$ is calm with respect to $p$ uniformly in $x$ at
$((\overline{x},\overline{p}),0),$ then $S$ is calm at $(\overline
{p},\overline{x})$ and%
\begin{equation}
\operatorname*{clm}S(\overline{p},\overline{x})\leq c^{-1}\widehat
{\operatorname*{clm}}_{p}H((\overline{x},\overline{p}),0). \label{mod_clmS}%
\end{equation}

(ii) If $H_{\overline{x}}$ is linearly pseudo-open with modulus $c>0$ at
$(\overline{p},0),$ then there exist $\gamma,\delta>0$ such that, for every
$p\in B(\overline{p},\gamma),$%
\begin{equation}
d(p,S^{-1}(\overline{x}))\leq c^{-1}d(0,H(\overline{x},p)\cap B(0,\delta)).
\label{subregS}%
\end{equation}

If, moreover, $H$ is calm with respect to $x$ uniformly in $p$ at
$((\overline{x},\overline{p}),0),$ then $S$ is metrically subregular at
$(\overline{p},\overline{x})$ and%
\begin{equation}
\operatorname*{subreg}S(\overline{p},\overline{x})\leq c^{-1}\widehat
{\operatorname*{clm}}_{x}H((\overline{x},\overline{p}),0). \label{mod_subregS}%
\end{equation}

\end{thm}

\noindent\textbf{Proof. }Observe first that is sufficient to prove the $(i)$
item, because the second item follows symmetrically, using $T:=S^{-1}$ instead
of $S,$ and taking into account Proposition \ref{link_at2}.

Let us prove the first part. We know that there exist $r,\varepsilon>0$ such
that, for every $\rho\in(0,\varepsilon)$ and every $x\in B(\overline{x},r)\cap
H_{\overline{p}}(B(0,c\rho)),$ one has $0\in H_{\overline{p}}(B(x,\rho)).$

Consider $\rho\in(0,\varepsilon),\alpha:=r,\beta:=c\rho,$ and take $x\in
B(\overline{x},\alpha).$ If $H(x,\overline{p})\cap B(0,\beta)=\emptyset$, the
relation (\ref{clmS}) automatically holds. Suppose next that $H(x,\overline
{p})\cap B(0,\beta)\not =\emptyset.$ If $0\in H(x,\overline{p})\cap
B(0,\beta),$ then, again, (\ref{clmS}) trivially holds. Consider now the case
$0\not \in H(x,\overline{p})\cap B(0,\beta).$ Then for every $\xi>0,$ one can
find $y_{\xi}\in H(x,\overline{p})\cap B(0,\beta)$ such that%
\[
\left\Vert y_{\xi}\right\Vert <d(0,H(x,\overline{p})\cap B(0,\beta))+\xi.
\]

Then%
\[
0\in B(y_{\xi},d(0,H(x,\overline{p})\cap B(0,\beta))+\xi).
\]

Because $d(0,H(x,\overline{p})\cap B(0,\beta))<\beta=c\rho,$ one can select
$\xi>0$ sufficiently small such that $d(0,H(x,\overline{p})\cap B(0,\beta
))+\xi<c\rho.$ Define $\rho_{0}:=c^{-1}(d(0,H(x,\overline{p})\cap
B(0,\beta))+\xi)<\rho<\varepsilon.$ Observe now that $y_{\xi}\in
H(x,\overline{p})\cap B(0,c\rho_{0}),$ which means that $x\in H_{\overline{p}%
}^{-1}(y_{\xi})\subset H_{\overline{p}}^{-1}(B(0,c\rho_{0})).$ Consequently,
$x\in B(\overline{x},r)\cap H_{\overline{p}}^{-1}(B(0,c\rho_{0})),$ so using
the assumption made we deduce that $0\in H_{\overline{p}}(B(x,\rho_{0})).$
Equivalently, there exists $\widetilde{x}\in B(x,\rho_{0})$ such that
$\widetilde{x}\in S(\overline{p}).$ In conclusion,%
\[
d(x,S(\overline{p}))\leq d(x,\widetilde{x})<\rho_{0}=c^{-1}d(0,H(x,\overline
{p})\cap B(0,\beta))+c^{-1}\xi.
\]

Making $\xi\rightarrow0,$ one obtains (\ref{clmS}).

Suppose next that $H$ is calm with respect to $p$ uniformly in $x$ at
$(\overline{x},\overline{p},0),$ so there exist $s,t,l>0$ such that
$ls<c\rho,$ and for every $(x,p)\in B(\overline{x},s)\times B(\overline
{p},s),$%
\[
e(H(x,p)\cap B(0,t),H(x,\overline{p}))\leq ld(p,\overline{p}).
\]

Consider $a:=\min\{\alpha,s\},$ and take $p\in B(\overline{p},s),x\in S(p)\cap
B(\overline{x},a).$ Then $0\in H(x,p)\cap B(0,t),$ so $d(0,H(x,\overline
{p}))\leq ld(p,\overline{p}).$ For every $\tau>0$ sufficiently small such that
$ls+\tau<c\rho,$ there is $y_{\tau}\in H(x,\overline{p})$ such that%
\[
\left\Vert y_{\tau}\right\Vert <ld(p,\overline{p})+\tau<c\rho.
\]

In conclusion, $y_{\tau}\in H(x,\overline{p})\cap B(0,c\rho),$ which means,
using (\ref{clmS}), that%
\[
d(x,S(\overline{p}))\leq c^{-1}d(0,H(x,\overline{p})\cap B(0,c\rho))\leq
c^{-1}\left\Vert y_{\tau}\right\Vert <c^{-1}ld(p,\overline{p})+c^{-1}\tau.
\]

Making $\tau\rightarrow0$ in the relation $d(x,S(\overline{p}))<c^{-1}%
ld(p,\overline{p})+c^{-1}\tau,$ and using the arbitrariness of $x\in S(p)\cap
B(\overline{x},a),$ one deduces the calmness of $S$ at $(\overline
{p},\overline{x}).$ Also, the relation between the associated moduli of
calmness easily follows.\hfill$\square$

\bigskip

Notice that Theorem \ref{main} (i) could be compared with \cite[Theorem
3.1]{CKY}, where the same conclusion is obtained under somehow stronger
assumptions in terms of coderivatives and using a closed-graph assumption for
$H.$ Moreover, our result could be deduced using even weaker concepts of
openness, but we preferred the actual form for consistency with results in the sequel.

\bigskip

The next examples emphasize the fact that in Theorem \ref{main} the converses
do not hold.

\begin{examp}
Consider the multifunctions $H:\mathbb{R\times R\rightrightarrows R},$ given
by%
\[
H(x,p):=\left\{
\begin{array}
[c]{ll}%
\{0\}, & \text{if }\left\vert x\right\vert \geq\left\vert p\right\vert \\
\left\{  \sqrt{\left\vert p\right\vert }\right\}  , & \text{if }\left\vert
x\right\vert <\left\vert p\right\vert .
\end{array}
\right.
\]
Then $S:\mathbb{R\rightrightarrows R}$ is given by $S(p)=\mathbb{R\setminus
(}-\left\vert p\right\vert ,\left\vert p\right\vert \mathbb{)}.$ Take
$(\overline{x},\overline{p})=(0,0).$ It is easy to prove that $S$ is
metrically subregular at\ $(0,0)$ with constant $1$ and calm at $(0,0)$ with
constant $1.$

On the other hand, $H_{\overline{x}}(p)=\left\{  \sqrt{\left\vert p\right\vert
}\right\}  $ for any $p$ and $H_{\overline{p}}(x)=\{0\}$ for any $x,$ whence
$H_{\overline{x}}$ is linearly pseudo-open at $(\overline{p},0)$ and
$H_{\overline{p}}$ is linearly pseudo-open at $(\overline{x},0).$

Finally, observe that $H$ is not calm with respect to $x$ uniformly in $p$ at
$((0,0),0)\ $and $H$ is not calm with respect to $p$ uniformly in $x$ at
$((0,0),0).$ If we suppose, by way of contradiction, that $H$ is calm with
respect to $x$ uniformly in $p$ at $((0,0),0)$ then there exist $L>0,$
$U\in\mathcal{V}(0),$ $V\in\mathcal{V}(0)$ and $W\in\mathcal{V}(0)$ such that
for every $x\in U$ and$\ p\in V$%
\[
e(H(x,p)\cap W,H(0,p))\leq Ld(x,0),
\]
hence we find $n_{0}\in\mathbb{N}$ such that $\sqrt{\frac{1}{n}}\leq L\frac
{1}{n}$ for every $n\geq n_{0}$, which is not true. Therefore, $H$ is not calm
with respect to $x$ uniformly in $p$ at $((0,0),0).$ Similarly, taking
$p=\frac{1}{n}$ and $x=\frac{1}{n^{2}}$ we deduce that $H$ is not calm with
respect to $p$ uniformly in $x$ at $((0,0),0).$
\end{examp}

\begin{examp}
Consider the multifunction $H:\mathbb{R\times R\rightrightarrows R},$ given by%
\[
H(x,p):=\left\{
\begin{array}
[c]{ll}%
\lbrack0,1], & \text{if }\left\vert x\right\vert \geq\left\vert p\right\vert
\\
(0,1], & \text{if }\left\vert x\right\vert <\left\vert p\right\vert .
\end{array}
\right.
\]

Then $S$ is the same as in the previous example, whence it is metrically
subregular at $(\overline{x},\overline{p})=(0,0)$ but $H_{\overline{x}}$ is
not linearly pseudo-open at $(0,0).$

Since $H_{0}(p)$ is $[0,1]$ if $\left\vert p\right\vert =0\ $and $(0,1]$ if
$\left\vert p\right\vert \neq0,$ $H_{0}$ is not linearly pseudo-open at
$(0,0)$ because if we suppose the opposite, taking $x=\frac{1}{n}$,
$y=\frac{1}{n^{2}}$ and $\rho=\frac{1}{Ln^{2}}$ we find $n_{0}\in\mathbb{N}$
such that $\frac{1}{n}<\frac{1}{Ln^{2}}$ for every $n\geq n_{0}$, which is not
true. In contrast, $H$ is calm with respect to $x$ uniformly in $p$ at
$((0,0),0).$
\end{examp}

\begin{rmk}
In both items of Theorem \ref{main}, one can replace the calmness assumptions
by somehow weaker (but more technical) conditions, as follows:

(i) Suppose that (\ref{clmS}) holds. If there exists $M>0$ such that for every
$x\in B(\overline{x},\alpha)$ and for every $p$ in a neighborhood of
$\overline{p}$ s.t. $d(p,\overline{p})<Md(0,H(x,\overline{p})\cap B(0,\beta))$
it follows that $0\notin H(x,p),$ then $S$ is calm at $(\overline{p}%
,\overline{x})$ with constant $c^{-1}M^{-1};$

(ii) Suppose that (\ref{subregS}) holds. If there exists $M>0$ such that for
every $p\in B(\overline{p},\gamma)$ and for every $x$ in a neighborhood of
$\overline{x}$ s.t. $d(x,\overline{x})<Md(0,H(\overline{x},p)\cap
B(0,\delta))$ it follows that $0\notin H(x,p),$ then $S$ is metrically
subregular at $(\overline{p},\overline{x})$ with constant $c^{-1}M^{-1}.$

Indeed, for (i), take $x\in B(\overline{x},\alpha)$ and $p$ close to
$\overline{p}$ s.t. $d(x,S(\overline{p}))>c^{-1}M^{-1}d(p,\overline{p}).$
Then
\[
c^{-1}d(0,H(x,\overline{p})\cap B(0,\beta))\geq d(x,S(\overline{p}%
))>c^{-1}M^{-1}d(p,\overline{p}),
\]
therefore $0\notin H(x,p),\ $i.e. $x\notin S(p).$ Taking into account Remark
\ref{char}, $S$ is calm with constant $c^{-1}M^{-1}$ at $(\overline
{p},\overline{x}).$ The second item follows symmetrically.

Observe that there are situations where this remark applies while Theorem
\ref{main} does not. This is shown in the example below.

Consider the multifunction $H:\mathbb{R\times R\rightrightarrows R},$ given by%
\[
H(x,p):=\left\{
\begin{array}
[c]{ll}%
\left\{  \frac{\left\vert x\right\vert }{\left\vert p\right\vert }\right\}
, & \text{if }p\neq0\text{ and }x\neq0,\\
\{0\}, & \text{if }p=0\text{ and }x\neq0,\\
\left\{  \left\vert p\right\vert \right\}  , & \text{if }x=0
\end{array}
\right.
\]
and take $\overline{x}=\overline{p}=0.$ Then $H_{\overline{x}}(p)=\left\{
\left\vert p\right\vert \right\}  $ is linearly pseudo-open with modulus $c=1$
at $(\overline{p},0)$ and there exist $\gamma=\delta=1$ such that for every
$p\in B(\overline{p},\gamma),$%
\[
d(p,S^{-1}(\overline{x}))=\left\vert p\right\vert =d(0,H(\overline{x},p)\cap
B(0,\delta)).
\]
For $M=1,$ if $d(x,\overline{x})<Md(0,H(\overline{x},p)\cap B(0,\delta))$ then
$\left\vert x\right\vert <\left\vert p\right\vert <1,$ whence $0\notin H(x,p)$
because $0\in H(x,p)$ if and only if $p=0.$ Then $S$ is metrically subregular
at $(0,0)$ according to the remark. For this case we cannot apply Theorem
\ref{main} because $H$ is not calm with respect to $x$ uniformly in $p$ at
$((0,0),0).$ To observe this, take $x=p=\frac{1}{n}$ and, arguing by
contradiction, one obtains $1-\frac{1}{n}\leq L\frac{1}{n}$ for every
$n\in\mathbb{N}$ large enough, which is not possible.
\end{rmk}

\subsection{Local sum-stability}

This subsection revisits the concept of local-sum stability introduced in
\cite[Section 4]{DurStr4}. Originally, this was used in relation with the
Aubin property of the sum-multifunction, while here we follow a similar
procedure, but for calmness. The notion itself reads as follows.

\begin{df}
\label{sum-sta}(\cite[Definition 4.2]{DurStr4}) Let $F:X\rightrightarrows Y,$
$G:X\rightrightarrows Y$ be two multifunctions and $(\overline{x}%
,{\overline{y},}\overline{z})\in X\times Y\times Y$ such that ${\overline
{y}\in F(\overline{x}),}$ ${\overline{z}\in G(\overline{x}).}$ We say that the
multifunction $(F,G)$ is locally sum-stable around $(\overline{x}%
,{\overline{y},}\overline{z})$ if for every $\varepsilon>0$ there exists
$\delta>0$ such that, for every $x\in B(\overline{x},\delta)$ and every
$w\in(F+G)(x)\cap B({\overline{y}+}\overline{z},\delta),$ there exist $y\in
F(x)\cap B({\overline{y},}\varepsilon)$ and $z\in G(x)\cap B(\overline
{z},\varepsilon)$ such that $w=y+z.$
\end{df}

\bigskip

Besides the initial results involving this notion (see \cite[Section
4]{DurStr4}), it was recently used and studied in relation with the metric
regularity of the sum of multifunctions in \cite{NTT2}.

\bigskip

We begin our analysis announced before by considering an example which shows
that the calmness property is not stable under summation (see, for more
details, \cite[Example 4.8]{DurStr4}, where a similar example is given to
prove that the Aubin property does not hold for the sum of multimappings).
Take the multifunctions $F,G:\mathbb{R\rightrightarrows R}$, given by%
\[
F(x):=\left\{
\begin{array}
[c]{ll}%
\lbrack0,1]\cup\{2\}, & \text{if }x\in\mathbb{R}\setminus\{0\}\\
\lbrack0,1], & \text{if }x=0
\end{array}
\right.
\]
and by $G(x):=[0,1]$ for every $x\in\mathbb{R}$, which are calm at $(0,1)$.
Then the multifunction $F+G:\mathbb{R}\rightrightarrows\mathbb{R}$, given by%
\[
(F+G)(x)=\left\{
\begin{array}
[c]{ll}%
\lbrack0,3], & \text{if }x\in\mathbb{R}\setminus\{0\}\\
\lbrack0,2], & \text{if }x=0,
\end{array}
\right.
\]

\noindent is not calm at $(0,2).$

\bigskip

The next lemma, whose proof is straightforward, shows that, as in the case of
Aubin property, the local sum-stability is the missing ingredient in order to
get the conservation of the calmness property at summation.

\begin{lm}
\label{clm_H}Let $F:X\rightrightarrows Y,$ $G:X\rightrightarrows Y$ be two
multifunctions. Suppose that $F$ is calm at $(\overline{x},{\overline{y}}%
)\in\operatorname{Gr}F,$ that $G$ is calm at $(\overline{x},{\overline{z}}%
)\in\operatorname{Gr}G,$ and that $(F,G)$ is locally sum-stable around
$(\overline{x},{\overline{y},}\overline{z}).$ Then the multifunction $F+G$ is
calm at $(\overline{x},{\overline{y}+}\overline{z}).$ Moreover, the following
relation holds true%
\begin{equation}
\operatorname*{clm}(F+G)(\overline{x},{\overline{y}+}\overline{z}%
)\leq\operatorname*{clm}F(\overline{x},{\overline{y}})+\operatorname*{clm}%
G(\overline{x},\overline{z}). \label{clm_sum}%
\end{equation}

\end{lm}

\noindent\textbf{Proof. }Use the calmness properties of $F$ and $G,$ to get
$\alpha,l,k>0$ such that, for every $x\in B(\overline{x},\alpha),$%
\begin{align}
e(F(x)\cap B({\overline{y},}\alpha),F(\overline{x}))  &  \leq ld(x,\overline
{x}),\label{clmF}\\
e(G(x)\cap B({\overline{z},}\alpha),G(\overline{x}))  &  \leq kd(x,\overline
{x}). \label{clmG}%
\end{align}

Using the locally sum-stability for $\varepsilon:=\alpha>0,$ one can find
$\delta\in(0,\alpha)$ such that, for every $x\in B(\overline{x},\delta)$ and
every $w\in(F+G)(x)\cap B({\overline{y}+}\overline{z},\delta),$ there exist
$y\in F(x)\cap B({\overline{y},}\alpha)$ and $z\in G(x)\cap B(\overline
{z},\alpha)$ such that $w=y+z.$ Consequently, using (\ref{clmF}) and
(\ref{clmG}),
\[
d(w,(F+G)(\overline{x}))\leq d(y,F(\overline{x}))+d(z,G(\overline{x}%
))\leq(l+k)d(x,\overline{x}).
\]
The relation (\ref{clm_sum}) easily follows.$\hfill\square$

\bigskip

Another remark is that, unsurprisingly, the calmness of the sum-multifunction
can be obtained in various situations, without the calmness of the component
multifunctions and in the absence of any kind of local sum-stability.

The next proposition, whose proof is omitted being again straightforward, uses
another sort of stability for the component mappings in proving the calmness
of the sum.

\begin{pr}
Let $F:X\rightrightarrows Y,$ $G:X\rightrightarrows Y$ be two multifunctions
and $(\overline{x},{\overline{y}})\in\operatorname{Gr}F,$ $(\overline
{x},\overline{z})\in\operatorname*{Gr}G$. Suppose that there exist
$L_{F},L_{G}>0$ and $U\in\mathcal{V}(\overline{x}),$ $V\in\mathcal{V}%
(\overline{y}+\overline{z})$ such that, for every $x\in U$ and $w\in V$ with
$(x,w)\in\operatorname*{Gr}(F+G),$ we can find $y\in F(x)$ with
$d(y,F(\overline{x}))\leq L_{F}d(x,\overline{x}),\ $and $z\in G(x)$ with
$d(y,F(\overline{x}))\leq L_{F}d(x,\overline{x}),\ $such that $w=y+z.$ Then
$F+G$ is $(L_{F}+L_{G})-$calm at $(\overline{x},{\overline{y}+}\overline{{z}%
}).$
\end{pr}

\bigskip

Remark that if $F$ is $L_{F}-$calm at $(\overline{x},\overline{y}),$ $G$ is
$L_{G}-$calm at $(\overline{x},\overline{z}),$ and $(F,G)$ is locally
sum-stable around $(\overline{x},\overline{y},\overline{z}),$ we are in the
setting of the previous proposition. Other situations are presented in the
following example.

\begin{examp}
1. Consider the multifunctions $F,G:\mathbb{R\rightrightarrows R},$ given by%
\[
F(x):=\left\{
\begin{array}
[c]{ll}%
\lbrack0,2], & \text{if }x\in\mathbb{R}\setminus\{0\}\\
\lbrack0,1], & \text{if }x=0
\end{array}
\right.
\]

\noindent and%
\[
G(x):=\left\{
\begin{array}
[c]{ll}%
\lbrack0,2], & \text{if }x\in\mathbb{R}\setminus\{0\}\\
\lbrack1,2], & \text{if }x=0,
\end{array}
\right.
\]

\noindent which are not calm at $(0,1)$. Then the multifunction
$F+G:\mathbb{R}\rightrightarrows\mathbb{R},$ given by%
\[
(F+G)(x)=\left\{
\begin{array}
[c]{ll}%
\lbrack0,4], & \text{if }x\in\mathbb{R}\setminus\{0\}\\
\lbrack1,3], & \text{if }x=0,
\end{array}
\right.
\]

\noindent is calm at $(0,2).$

2. Consider the multifunctions $F,G:\mathbb{[-}1,1\mathbb{]\rightrightarrows
R},$ given by%
\[
F(x):=\left\{
\begin{array}
[c]{ll}%
\lbrack0,x+1], & \text{if }x\in\lbrack-1,0)\\
\{0\}\cup\left[  \dfrac{1}{2},1\right]  , & \text{if }x=0\\
\lbrack0,1-x], & \text{if }x\in(0,1],
\end{array}
\right.
\]

\noindent which is not calm at $(0,0),$ and%
\[
G(x):=\left\{
\begin{array}
[c]{ll}%
\{-1-x,0\}, & \text{if }x\in\lbrack-1,0)\\
\lbrack-1,0], & \text{if }x=0\\
\{-1+x,0\}, & \text{if }x\in(0,1],
\end{array}
\right.
\]

\noindent which is calm at $(0,0)$. Then the multifunction $F+G:\mathbb{R}%
\rightrightarrows\mathbb{R},$ given by%
\[
(F+G)(x)=\left\{
\begin{array}
[c]{ll}%
\lbrack-x-1,x+1], & \text{if }x\in\lbrack-1,0)\\
\lbrack-1,1], & \text{if }x=0\\
\lbrack x-1,1-x], & \text{if }x\in(0,1]
\end{array}
\right.
\]

\noindent is calm at $(0,0).$ Remark that $(F,G)$ is not locally sum-stable
around $(0,0,0)$.
\end{examp}

\subsection{Applications to variational systems}

This subsection plays a leading role in this work, since here we put together
all the facts collected by now and we use them in order to study variational systems.

\bigskip

To begin, we adapt the definition of local-sum stability to the parametric case.

\begin{df}
\label{sum-sta_param}(\cite[Definition 4.9]{DurStr4}) Let $F:X\times
P\rightrightarrows Y,$ $G:X\rightrightarrows Y$ be two multifunctions and
$(\overline{x},\overline{p},{\overline{y},}\overline{z})\in X\times P\times
Y\times Y$ such that ${\overline{y}\in F(\overline{x},\overline{p}),}$
${\overline{z}\in G(\overline{x}).}$ We say that the multifunction $(F,G)$ is
locally sum-stable around $((\overline{x},\overline{p}),{\overline{y}%
,}\overline{z})$ if for every $\varepsilon>0$ there exists $\delta>0$ such
that, for every $(x,p)\in B(\overline{x},\delta)\times B(\overline{p},\delta)$
and every $w\in(F_{p}+G)(x)\cap B({\overline{y}+}\overline{z},\delta),$ there
exist $y\in F_{p}(x)\cap B({\overline{y},}\varepsilon)$ and $z\in G(x)\cap
B(\overline{z},\varepsilon)$ such that $w=y+z.$
\end{df}

Also, Lemma \ref{clm_H} has the following variant in the parametric case.

\begin{lm}
\label{clm_H_par}Let $F:X\times P\rightrightarrows Y,$ $G:X\rightrightarrows
Y$ be two multifunctions. Suppose that $F$ is calm with respect to $x$
uniformly in $p$ at $((\overline{x},\overline{p}),{\overline{y}}%
)\in\operatorname{Gr}F,$ that $G$ is calm at $(\overline{x},{\overline{z}}%
)\in\operatorname{Gr}G$ and that $(F,G)$ is locally sum-stable around
$((\overline{x},\overline{p}),{\overline{y},}\overline{z}).$ Then the
multifunction $H:X\times P\rightrightarrows Y$ given by $H(x,p):=F(x,p)+G(x)$
is calm with respect to $x$ uniformly in $p$ at $((\overline{x},\overline
{p}),{\overline{y}+}\overline{z}).$ Moreover, the following relation holds
true%
\begin{equation}
\widehat{\operatorname*{clm}}_{x}H((\overline{x},\overline{p}),{\overline{y}%
+}\overline{z})\leq\widehat{\operatorname*{clm}}_{x}F((\overline{x}%
,\overline{p}),{\overline{y}})+\operatorname*{clm}G(\overline{x},\overline
{z}). \label{clm_sum_par}%
\end{equation}

\end{lm}

\noindent\textbf{Proof.} Adapt the line of the proof of Lemma \ref{clm_H}%
.\hfill$\square$

\bigskip

In the next theorems we get at-point regularity results for $S$ defined by
(\ref{S}), where $H$ takes the form
\begin{equation}
H(x,p):=F(x,p)+G(x). \label{H}%
\end{equation}
The case of around-point regularity was considered in \cite{ArtMord2009},
\cite{DurStr4}, while particular situations for at-point regularity are
studied in \cite{ArtMord2009}, \cite{ArtMord2010}.

\begin{thm}
\label{msubreg_sol}Let $X,Y,P$ be Banach spaces, $F:X\times P\rightrightarrows
Y,$ $G:X\rightrightarrows Y$ be two set-valued maps and $(\overline
{x},\overline{p},\overline{y})\in X\times P\times Y$ such that $\overline
{y}\in F(\overline{x},\overline{p})$ and $-\overline{y}\in G(\overline{x})$.
Suppose that the following assumptions are satisfied:

(i) $(F,G)$ is locally sum-stable around $((\overline{x},\overline
{p}),{\overline{y},}-{\overline{y}});$

(ii) $F$ is calm with respect to $x$ uniformly in $p$ at $((\overline
{x},\overline{p}),\overline{y});$

(iii) $F_{\overline{x}}$ is metrically regular around $(\overline{p}%
,\overline{y});$

(iv) $G$ is calm at $(\overline{x},-{\overline{y}}).$

Then $S$ is metrically subregular at $(\overline{p},\overline{x}).$ Moreover,
the next relation holds%
\begin{equation}
\operatorname*{subreg}S(\overline{p},\overline{x})\leq\operatorname*{reg}%
F_{\overline{x}}(\overline{p},\overline{y})\cdot\lbrack\widehat
{\operatorname*{clm}}_{x}F((\overline{x},\overline{p}),\overline
{y})+\operatorname*{clm}G(\overline{x},-{\overline{y}})]. \label{rS}%
\end{equation}

\end{thm}

\noindent\textbf{Proof. }Using Lemma \ref{clm_H_par}, we know that $H$ given
by (\ref{H}) is calm with respect to $x$ uniformly in $p$ at $((\overline
{x},\overline{p}),{0})$ and the relation (\ref{clm_sum_par}) holds for
$\overline{z}:=-{\overline{y}.}$

Using (iii), which is equivalent to the linear openness of $F_{\overline{x}}$
around $(\overline{p},\overline{y}),$ one can find $\alpha,L>0$ such that, for
every $(p,y)\in\operatorname*{Gr}F_{\overline{x}}\cap\lbrack B(\overline
{p},\alpha)\times B(\overline{y},\alpha)]$ and every $\rho\in(0,\alpha),$%
\[
B(y,L\rho)\subset F_{\overline{x}}(B(p,\rho)).
\]

Also, from (i), for $\alpha$ found before, there is $\delta\in(0,\alpha)$ such
that, for every $(x,p)\in B(\overline{x},\delta)\times B(\overline{p},\delta)$
and every $w\in H(x,p)\cap B(0,\delta),$ one can find $y\in F(x,p)\cap
B(\overline{y},\alpha)$ and $z\in G(x)\cap B(-\overline{y},\alpha)$ such that
$w=y+z.$

Fix now $\tau<\min\{\alpha,L^{-1}\delta\},$ and take $\rho\in(0,\tau),$ $p\in
B(\overline{p},\delta)\cap H_{\overline{x}}^{-1}(B(0,L\rho)).$ Then there
exists $w\in B(0,L\rho)$ such that $w\in H(\overline{x},p).$ Then $w\in
B(0,\delta),$ so using the local sum-stability of $(F,G),$ one can find $y\in
F(\overline{x},p)\cap B(\overline{y},\alpha)$ and $z\in G(\overline{x})\cap
B(-\overline{y},\alpha)$ such that $w=y+z.$ Then%
\[
0\in B(w,L\rho)=B(y,L\rho)+z\subset F_{\overline{x}}(B(p,\rho))+z\subset
H_{\overline{x}}(B(p,\rho)).
\]

As consequence, $H_{\overline{x}}$ is linearly pseudo-open at $(\overline
{p},0),$ with modulus $(\operatorname*{reg}F_{\overline{x}}(\overline
{p},\overline{y}))^{-1}.$ The conclusion now follows from the second part of
Theorem \ref{main}.$\hfill\square$

\bigskip

A natural question which arises when one looks at Theorem \ref{msubreg_sol} is
if the assumption (iii) cannot be weakened, supposing for example just the
metric subregularity of $F_{\overline{x}}$ at $(\overline{p},\overline{y}).$
The next example clarifies this aspect, showing that if one replaces the
metric regularity of $F_{\overline{x}}$ with its metric subregularity, the
conclusion of Theorem \ref{msubreg_sol} is not satisfied in general.

\begin{examp}
Consider the multifunctions $F:\mathbb{R}^{2}\rightrightarrows\mathbb{R}^{2}$
and $G:\mathbb{R}\rightrightarrows\mathbb{R}^{2}$ given by%
\begin{align*}
F(x,p)  &  :=\{(0,p)\}\cup\left\{  \left(  \frac{1}{n^{2}},p-\frac{1}{n^{2}%
}\right)  \mid n\in\mathbb{N}\setminus\{0\}\right\}  \text{ and}\\
G(x)  &  :=\{(x,0)\}\cup\left\{  \left(  x+\frac{1}{n^{3}},\frac{1}{n}\right)
\mid n\in\mathbb{N}\setminus\{0\}\right\}  .
\end{align*}

\noindent Also, fix $\overline{x}:=0,$ $\overline{p}:=0,$ $\overline
{y}:=(0,0).$ Then $H:\mathbb{R}^{2}\rightrightarrows\mathbb{R}^{2}$ is given
by%
\begin{align*}
H(x,p)  &  =\{(x,p)\}\cup\left\{  \left(  x+\frac{1}{n^{2}},p-\frac{1}{n^{2}%
}\right)  \mid n\in\mathbb{N}\setminus\{0\}\right\}  \cup\\
&  \left\{  \left(  x+\frac{1}{n^{3}},p+\frac{1}{n}\right)  \mid
n\in\mathbb{N}\setminus\{0\}\right\}  \cup\\
&  \left\{  \left(  x+\frac{1}{n^{3}}+\frac{1}{m^{2}},p-\frac{1}{m^{2}}%
+\frac{1}{n}\right)  \mid n,m\in\mathbb{N}\setminus\{0\}\right\}  .
\end{align*}

Then,%
\[
S^{-1}(x)=\left\{
\begin{array}
[c]{l}%
0,\text{ if }x=0\\
\dfrac{1}{n^{2}},\text{ if }x=-\dfrac{1}{n^{2}},n\in\mathbb{N}\setminus\{0\}\\
-\dfrac{1}{n},\text{ if }x=-\dfrac{1}{n^{3}},n\in\mathbb{N}\setminus\{0\}\\
\dfrac{1}{m^{2}}-\dfrac{1}{n},\text{ if }x=-\dfrac{1}{n^{3}}-\dfrac{1}{m^{2}%
},m,n\in\mathbb{N}\setminus\{0\}\\
\emptyset,\text{ otherwise.}%
\end{array}
\right.
\]

Let us remark that $S^{-1}(0)=\{0\}.$ Also, one can prove that $S^{-1}$ is not
calm at $(0,0).$ If $S^{-1}$ would be calm at $(0,0),$ it should exist
$l,\alpha>0$ such that, for every $x\in B(0,\alpha),$%
\[
S^{-1}(x)\cap B(0,\alpha)\subset S^{-1}(0)+l\left\vert x\right\vert [-1,1].
\]

\noindent Take $x=-\dfrac{1}{n^{3}}-\dfrac{1}{m^{2}},$ with $m,n\in
\mathbb{N}\setminus\{0\},$ such that $\dfrac{1}{n^{3}}+\dfrac{1}{m^{2}}%
<\alpha$ and $\left\vert \dfrac{1}{m^{2}}-\dfrac{1}{n}\right\vert <\alpha.$
One should have%
\[
\left\vert \frac{1}{m^{2}}-\frac{1}{n}\right\vert <l\left(  \frac{1}{n^{3}%
}+\frac{1}{m^{2}}\right)  ,
\]

\noindent for every $m,n$ sufficiently large. But for $m=n,$ one should have
that%
\begin{align*}
\frac{1}{n}-\frac{1}{n^{2}}  &  <l\left(  \frac{1}{n^{3}}+\frac{1}{n^{2}%
}\right)  ,\text{ or}\\
\frac{(n-1)n}{n+1}  &  <l
\end{align*}

\noindent for every $n$ sufficiently large, which is absurd. In conclusion,
$S^{-1}$ is not calm at $(0,0),$ or $S$ is not metrically subregular at
$(0,0).$

Let us prove now that $F_{0}$ is metrically subregular at $(0,(0,0)),$ but $F$
is not metrically regular around $(0,(0,0)).$

As one can see,%
\[
F_{0}^{-1}(u,v)=\left\{
\begin{array}
[c]{l}%
p,\text{ if }(u,v)=(0,p)\text{ or }(u,v)=\left(  \dfrac{1}{n^{2}},p-\dfrac
{1}{n^{2}}\right)  ,\text{ }n\in\mathbb{N}\setminus\{0\}\\
\emptyset,\text{ otherwise.}%
\end{array}
\right.
\]

\noindent Then $F_{0}^{-1}(0,0)=\{0\}.$ Also, let's prove that there exists
$\beta>0$ such that, for every $(u,v)\in B((0,0),\beta),$%
\begin{equation}
F_{0}^{-1}(u,v)\cap B(0,\beta)\subset F_{0}^{-1}(0,0)+\left\Vert
(u,v)\right\Vert \mathbb{[-}1,1]. \label{F0}%
\end{equation}

\noindent Indeed, take arbitrary $\beta>0$ and $(u,v)\in B((0,0),\beta)$ such
that $F_{0}^{-1}(u,v)\not =\emptyset.$ Then $(u,v)=(0,p),$ in which case
relation (\ref{F0}) reduces to%
\[
\left\vert p\right\vert \leq\left\Vert (0,p)\right\Vert ,
\]

\noindent or $(u,v)=\left(  \dfrac{1}{n^{2}},p-\dfrac{1}{n^{2}}\right)  ,$ in
which case relation (\ref{F0}) becomes%
\[
\left\vert p\right\vert \leq\frac{1}{n^{2}}+\left\vert p-\frac{1}{n^{2}%
}\right\vert ,
\]

\noindent which is obviously true from the triangle inequality.

But, as one can easily see, $F_{0}^{-1}(u^{\prime},v^{\prime})$ can be empty
for $(u^{\prime},v^{\prime})$ arbitrarily close to $(0,0),$ hence the relation%
\[
F_{0}^{-1}(u,v)\cap B(0,\beta)\subset F_{0}^{-1}(u^{\prime},v^{\prime
})+m\left\Vert (u,v)-(u^{\prime},v^{\prime})\right\Vert \mathbb{[-}1,1]
\]

\noindent cannot be true for every $(u^{\prime},v^{\prime})$ in a neighborhood
of $(0,0).$ As consequence, $F_{0}^{-1}$ does not have the Aubin property
$((0,0),0),$ or $F_{0}$ is not metrically regular around $(0,(0,0)).$

Being constant with respect to $x,$ $F$ is obviously calm with respect to $x,$
uniformly in $p$ at $((0,0),(0,0))$ with modulus $1.$ Also, one can easily see
that $G$ is calm at $(0,(0,0))$ with modulus $1.$

Let us now prove that $(F,G)$ is locally sum-stable around
$((0,0),(0,0),(0,0)).$ Indeed, take arbitrary $\varepsilon>0.$ Pick
$\delta<\min\left\{  \dfrac{2}{7}\varepsilon,\dfrac{18}{343}\varepsilon
^{3}\right\}  ,$ and take $x\in\left(  -\dfrac{\delta}{2},\dfrac{\delta}%
{2}\right)  ,$ $p\in\left(  -\dfrac{\delta}{2},\dfrac{\delta}{2}\right)  $ and
$w\in H(x,p)\cap B((0,0),\delta).$ We have four possibilities:

1. If $w=(x,p),$ then the conclusion easily follows.

2. If $w=\left(  x+\dfrac{1}{n^{2}},p-\dfrac{1}{n^{2}}\right)  ,$ with
$\left\vert x+\dfrac{1}{n^{2}}\right\vert +\left\vert p-\dfrac{1}{n^{2}%
}\right\vert <\delta,$ then $\dfrac{1}{n^{2}}<\dfrac{3}{2}\delta,$ because
otherwise
\[
\left\vert x+\frac{1}{n^{2}}\right\vert \geq\frac{1}{n^{2}}-\left\vert
x\right\vert \geq\frac{3}{2}\delta-\frac{1}{2}\delta=\delta,
\]

\noindent which is absurd. Then $w$ can be written as $\left(  \dfrac{1}%
{n^{2}},p-\dfrac{1}{n^{2}}\right)  +(x,0),$ with $\left(  \dfrac{1}{n^{2}%
},p-\dfrac{1}{n^{2}}\right)  \in F(x,p),$ $(x,0)\in G(x),$ and%
\begin{align*}
\left\Vert \left(  \frac{1}{n^{2}},p-\frac{1}{n^{2}}\right)  \right\Vert  &
=\frac{1}{n^{2}}+\left\vert p-\frac{1}{n^{2}}\right\vert <\frac{5}{2}%
\delta<\varepsilon,\\
\left\Vert (x,0)\right\Vert  &  =\left\vert x\right\vert <\frac{1}{2}%
\delta<\varepsilon.
\end{align*}

3. If $w=\left(  x+\dfrac{1}{n^{3}},p+\dfrac{1}{n}\right)  ,$ with $\left\vert
x+\dfrac{1}{n^{3}}\right\vert +\left\vert p-\dfrac{1}{n}\right\vert <\delta,$
then $\dfrac{1}{n}<\dfrac{3}{2}\delta,$ because, again, in the opposite
situation one should have
\[
\left\vert p-\frac{1}{n}\right\vert \geq\frac{1}{n}-\left\vert p\right\vert
\geq\frac{3}{2}\delta-\frac{1}{2}\delta=\delta,
\]

\noindent which is absurd. Then $w$ can be written as $\left(  0,p\right)
+\left(  x+\dfrac{1}{n^{3}},\dfrac{1}{n}\right)  ,$ with $\left(  0,p\right)
\in F(x,p),$ $\left(  x+\dfrac{1}{n^{3}},\dfrac{1}{n}\right)  \in G(x),$ and%
\begin{align*}
\left\Vert \left(  0,p\right)  \right\Vert  &  =\left\vert p\right\vert
<\frac{1}{2}\delta<\varepsilon,\\
\left\Vert \left(  x+\frac{1}{n^{3}},\frac{1}{n}\right)  \right\Vert  &
=\left\vert x+\frac{1}{n^{3}}\right\vert +\frac{1}{n}<\frac{5}{2}%
\delta<\varepsilon.
\end{align*}

4. Finally, if $w=\left(  x+\dfrac{1}{n^{3}}+\dfrac{1}{m^{2}},p-\dfrac
{1}{m^{2}}+\dfrac{1}{n}\right)  ,$ with $\left\vert x+\dfrac{1}{n^{3}}%
+\dfrac{1}{m^{2}}\right\vert +\left\vert p-\dfrac{1}{m^{2}}+\dfrac{1}%
{n}\right\vert <\delta,$ then, as above, one can prove that $\left\vert
\dfrac{1}{n^{3}}+\dfrac{1}{m^{2}}\right\vert <\dfrac{3}{2}\delta.$ But this
means that $\dfrac{1}{n}<\sqrt[3]{\dfrac{3}{2}\delta}$ and $\dfrac{1}{m}%
<\sqrt{\dfrac{3}{2}\delta}.$ Then $w$ can be written as $\left(  \dfrac
{1}{m^{2}},p-\dfrac{1}{m^{2}}\right)  +\left(  x+\dfrac{1}{n^{3}},\dfrac{1}%
{n}\right)  ,$ with $\left(  \dfrac{1}{m^{2}},p-\dfrac{1}{m^{2}}\right)  \in
F(x,p),$ $\left(  x+\dfrac{1}{n^{3}},\dfrac{1}{n}\right)  \in G(x),$ and%
\begin{align*}
\left\Vert \left(  \frac{1}{m^{2}},p-\frac{1}{m^{2}}\right)  \right\Vert  &
=\frac{1}{m^{2}}+\left\vert p-\frac{1}{m^{2}}\right\vert \leq\frac{2}{m^{2}%
}+\left\vert p\right\vert <3\delta+\frac{\delta}{2}=\frac{7}{2}\delta
<\varepsilon,\\
\left\Vert \left(  x+\frac{1}{n^{3}},\frac{1}{n}\right)  \right\Vert  &
=\left\vert x+\frac{1}{n^{3}}\right\vert +\frac{1}{n}\leq\left\vert
x\right\vert +\frac{1}{n^{3}}+\frac{1}{n}<\frac{\delta}{2}+\frac{3}{2}%
\delta+\sqrt[3]{\frac{3}{2}\delta}<\varepsilon.
\end{align*}

\noindent In conclusion, $(F,G)$ is locally sum-stable around
$((0,0),(0,0),(0,0)).$

Let us finally prove that, for every $L>0,$ there exist $\rho$ arbitrary
small, $p$ arbitrary small, and $w\in H_{0}(p)\cap D((0,0),\rho)$ such that
$w$ cannot be written as $y+z,$ with $y\in F_{0}(p)\cap D((0,0),L\rho)$ and
$z\in G(0).$ Indeed, take $\rho=\dfrac{1}{n^{3}}+\dfrac{2}{n^{2}},$
$p=-\dfrac{1}{n}$ arbitrary small, and $w=\left(  \dfrac{1}{n^{3}}+\dfrac
{1}{n^{2}},-\dfrac{1}{n^{2}}\right)  \in H_{0}(p)$ such that $\left\Vert
w\right\Vert =\rho.$ This $w$ can only be obtained as $\left(  \dfrac{1}%
{n^{2}},-\dfrac{1}{n}-\dfrac{1}{n^{2}}\right)  +\left(  \dfrac{1}{n^{3}%
},\dfrac{1}{n}\right)  ,$ with $\left(  \dfrac{1}{n^{2}},-\dfrac{1}{n}%
-\dfrac{1}{n^{2}}\right)  \in F_{0}(p),$ and $\left(  \dfrac{1}{n^{3}}%
,\dfrac{1}{n}\right)  \in G(0).$ Let us prove this assertion. Suppose there
exist $k,m\in\mathbb{N}\setminus\{0\}$ such that%
\[
\left\{
\begin{array}
[c]{c}%
\dfrac{1}{k^{3}}+\dfrac{1}{m^{2}}=\dfrac{1}{n^{3}}+\dfrac{1}{n^{2}}\\
-\dfrac{1}{m^{2}}+\dfrac{1}{k}-\dfrac{1}{n}=-\dfrac{1}{n^{2}}%
\end{array}
\right.  .
\]

\noindent This means that $\dfrac{1}{k^{3}}+\dfrac{1}{k}=\dfrac{1}{n^{3}%
}+\dfrac{1}{n}.$ As the function $k\mapsto\dfrac{1}{k^{3}}+\dfrac{1}{k}$ is
strictly decreasing, the only solution of this equation is $k=n.$ But this
shows also that $m=n.$

Now, suppose that there is $L>0$ such that $\left\Vert \left(  \dfrac{1}%
{n^{2}},-\dfrac{1}{n}-\dfrac{1}{n^{2}}\right)  \right\Vert <L\rho=L\left(
\dfrac{1}{n^{3}}+\dfrac{2}{n^{2}}\right)  ,$ for every $n$ sufficiently large.
This means that%
\begin{align*}
\frac{2}{n^{2}}+\frac{1}{n}  &  <L\left(  \frac{1}{n^{3}}+\frac{2}{n^{2}%
}\right)  ,\text{ or}\\
\frac{(n+2)n}{2n+1}  &  <L,
\end{align*}
for every $n$ sufficiently large, which is absurd.$\hfill\square$
\end{examp}

\bigskip

Notice that the phenomenon described at the final of the previous example
(which in fact generated it) seems to be the reason for which only the
subregularity condition for $F_{\overline{x}}$ is not sufficient in getting
the conclusion of Theorem \ref{msubreg_sol}. In fact, the problem was the
impossibility to obtain a (linear) correspondence needed to link the closeness
between $w$ and $\overline{w}$ with the one between $y$ and $\overline{y}.$ A
possible solution could be to introduce a notion of partial linear
sum-stability (in which to ask for this linear correspondence, but to maintain
only $y$ close to $\overline{y}$). We preferred to avoid this approach for
clarity. The next example, which seems to be more simple than the one given
before (at least in verifying the conditions for the involved objects) mainly
addresses the same questions, but the phenomenon described before is less
visible. For this reason, we keep them both.

\bigskip

\begin{examp}
Consider the multifunctions $F:\mathbb{R}^{2}\rightrightarrows\mathbb{R}$ and
$G:\mathbb{R}\rightrightarrows\mathbb{R}$ given by
\begin{align*}
F(x,p)  &  :=\left\{
\begin{array}
[c]{ll}%
\lbrack\left\vert p\right\vert ,+\infty)\cap%
\mathbb{Q}
, & \text{if }x=0\\
\lbrack\left\vert p\right\vert ,+\infty), & \text{if }x\text{ }\neq0,
\end{array}
\right. \\
G(x)  &  :=\left\{
\begin{array}
[c]{ll}%
\lbrack-1,0]\setminus%
\mathbb{Q}
\cup\{0\}, & \text{if }x=0\\
\lbrack-1,0], & \text{if }x\text{ }\neq0
\end{array}
\right.
\end{align*}
and take $\overline{x}=\overline{p}=\overline{y}=\overline{z}:=0.$ Then
$F\ $and $G$ here share all the properties from the previous example and,
again, $S$ is not metrically subregular at $(0,0).$
\end{examp}

The next theorem, previously given in \cite{DurStr4}, is a sort of
Lyusternik-Graves type result, and has the role to precisely specify the
constants involved in the openness property. This will be important in the
development of other subsequent results.

\begin{thm}
\label{main_const}(\cite[Theorem 3.3]{DurStr4}) Let $X,Y$ be Banach spaces,
$F_{1}:X\rightrightarrows Y$ and $F_{2}:Y\rightrightarrows X$ be two
multifunctions and $(\overline{x},{\overline{y}},\overline{z})\in X\times
Y\times Y$ such that $(\overline{x},{\overline{y}})\in\operatorname{Gr}F_{1}$
and $(\overline{z},\overline{x})\in\operatorname{Gr}F_{2}.$ Suppose that the
following assumptions are satisfied:

(i) $\operatorname{Gr}F_{1}$ is locally closed around $(\overline
{x},{\overline{y}}),$ so there exist $\alpha_{1},\beta_{1}>0$ such that
$\operatorname{Gr}F_{1}\cap\lbrack D(\overline{x},\alpha_{1})\times
D({\overline{y}},\beta_{1})]$ is closed;

(ii) $\operatorname{Gr}F_{2}$ is locally closed around $(\overline
{z},\overline{x}),$ so there exist $\alpha_{2},\beta_{2}>0$ such that
$\operatorname{Gr}F_{2}\cap\lbrack D(\overline{z},\beta_{2})\times
D(\overline{x},\alpha_{2})]$ is closed;

(iii) there exist $L,r_{1},s_{1}>0$ such that, for every $(x^{\prime
},y^{\prime})\in\operatorname{Gr}F_{1}\cap\lbrack B(\overline{x},r_{1})\times
B({\overline{y}},s_{1})],$ $F_{1}$ is $L-$open at $(x^{\prime},y^{\prime});$

(iv) there exist $M,r_{2},s_{2}>0$ such that, for every $(v^{\prime}%
,u^{\prime})\in\operatorname{Gr}F_{2}\cap\lbrack B(\overline{z},s_{2})\times
B(\overline{x},r_{2})],$ $F_{2}$ is $M-$open at $(v^{\prime},u^{\prime});$

(v) $LM>1.$

Then for every $\rho\in(0,\varepsilon),$ where $\varepsilon:=\min\{\alpha
_{1},\alpha_{2},L^{-1}\beta_{1},M\beta_{2},r_{1},r_{2},L^{-1}s_{1},Ms_{2}\},$%
\[
B({\overline{y}}-\overline{z},(L-M^{-1})\rho)\subset(F_{1}-F_{2}%
^{-1})(B(\overline{x},\rho)).
\]

\end{thm}

We are able to present now a fixed-point assertion, given in the parametric
form, which follows the path opened by Arutyunov in a series of papers
(\cite{Arut2007}--\cite{AZ}), and after that continued by Donchev and
Frankowska (\cite{DonFra2010}, \cite{DonFra2011}), Ioffe (\cite{Ioffe2010b}),
and Durea and Strugariu (\cite{DurStr5}, \cite{DurStr6}), where links between
fixed-point theorems and Lyusternik-Graves type results are provided. The
inequality from the conclusion of Theorem \ref{fixp} could be formulated in
some different ways, as is done in \cite[Theorem 7]{DonFra2010}, \cite[Theorem
4.4]{DurStr5}. We prefer to obtain it just in the form we need in the sequel.
For this, consider the multifunctions $\Phi:X\times P\rightrightarrows Y,$
$\Psi:X\rightrightarrows Y,$ and take the implicit multifunction
$S:P\rightrightarrows X$ as
\begin{align*}
S(p)  &  :=(\Phi(\cdot,p)-\Psi)^{-1}(0)=\{x\in X\mid0\in\Phi(x,p)-\Psi(x)\}\\
&  =\{x\in X\mid\Phi(x,p)\cap\Psi(x)\neq\emptyset\}=\operatorname*{Fix}%
(\Phi(\cdot,p)^{-1}\Psi).
\end{align*}

Although the proof of the next result has some common points to the one of the
first part of Theorem \ref{main} (i), it involves some more technicalities.
For this reason, we present it in full extent.

\begin{thm}
\label{fixp}Let $X,Y$ be Banach spaces, $P$ be a metric space, $\Phi:X\times
P\rightrightarrows Y$ and $\Psi:X\rightrightarrows Y$ be multifunctions and
$(\overline{x},\overline{p},\overline{y})\in X\times P\times Y$ such that
$((\overline{x},\overline{p}),\overline{y})\in\operatorname{Gr}\Phi$ and
$(\overline{x},\overline{y})\in\operatorname{Gr}\Psi.$ Suppose that the
following assumptions are satisfied:

(i) $\operatorname*{Gr}\Phi_{p}$ is locally closed around $(\overline
{x},\overline{y})$ uniformly for $p$ in a neighborhood of $\overline{p};$

(ii) $\operatorname*{Gr}\Psi$ is locally closed around $(\overline
{x},\overline{y});$

(iii) $(\Phi,-\Psi)$ is locally sum-stable around $((\overline{x},\overline
{p}),\overline{y},-\overline{y})$;

(iv) $\Phi$ has the Aubin property with respect to $x$ uniformly in $p$ around
$((\overline{x},\overline{p}),\overline{y})$ with constant $l>0;$

(v) $\Psi$ is metrically regular with constant $m>0$ around $(\overline
{x},\overline{y});$

(vi) $lm<1.$

Then there exist $\alpha,\beta>0$ such that for any $(x,p)\in B(\overline
{x},\alpha)\times B(\overline{p},\alpha)$ one has that
\begin{equation}
d(x,S(p))\leq(m^{-1}-l)^{-1}d(0,[\Phi(x,p)-\Psi(x)]\cap B(0,\beta)).
\label{diffix}%
\end{equation}

\end{thm}

\noindent\textbf{Proof.} Using (i) and (ii), one can pick $\gamma>0$ such
that, for every $p\in B(\overline{p},\gamma),$ $\operatorname*{Gr}\Phi_{p}%
\cap\lbrack D(\overline{x},\gamma)\times D(\overline{y},\gamma)]$ is closed
and $\operatorname*{Gr}\Psi\cap\lbrack D(\overline{x},\gamma)\times
D(\overline{y},\gamma)]$ is closed.

Also, by (iv), there exist $r\in(0,\gamma)$\ such that, for every $p\in
B(\overline{p},r),$ and every $x,x^{\prime}\in B(\overline{x},r),$%
\begin{equation}
e(\Phi_{p}(x)\cap B(\overline{y},r),\Phi_{p}(x^{\prime}))\leq ld(x,x^{\prime
}). \label{Aub_phip}%
\end{equation}

Take now $x\in B(\overline{x},2^{-1}r)$ and $y\in\Phi_{p}(x)\cap
B(\overline{y},r).$ Then for every $x^{\prime}\in B(x,2^{-1}r),$ we get by
(\ref{Aub_phip}) that%
\[
d(y,\Phi_{p}(x^{\prime}))\leq ld(x,x^{\prime}),
\]

\noindent which proves that for every $p\in B(\overline{p},r),$ $\Phi_{p}$ is
$l-$pseudocalm at every $(x,y)\in\operatorname*{Gr}\Phi_{p}\cap\lbrack
B(\overline{x},2^{-1}r)\times B(\overline{y},r)].$ In virtue of Proposition
\ref{link_at}, we deduce that that for every $p\in B(\overline{p},r),$
$\Phi_{p}^{-1}$ is $l^{-1}-$open at every $(y,x)\in\operatorname*{Gr}\Phi
_{p}^{-1}\cap\lbrack B(\overline{y},r)\times B(\overline{x},2^{-1}r)].$

From (v), we know that there exist $t>0$ such that, for every $(u,v)\in
\operatorname{Gr}\Psi\cap\lbrack B(\overline{x},t)\times B(\overline{y},t)],$
$\Psi$ is metrically hemiregular at $(u,v)$ with constant $m,$ hence is open
at linear rate $m^{-1}$ at $(u,v).$

From the local sum-stability of $(\Phi,-\Psi)$ around $((\overline
{x},\overline{p}),\overline{y},-\overline{y}),$ using $\min\{2^{-1}%
r,2^{-1}t\}$ instead of $\varepsilon,$ one can find $\delta>0$ such that, for
every $(x,p)\in B(\overline{x},\delta)\times B(\overline{p},\delta)$ and every
$w\in(\Phi_{p}-\Psi)(x)\cap B({0},\delta),$ there exist $y\in\Phi_{p}(x)\cap
B({\overline{y},}\varepsilon)$ and $z\in\Psi(x)\cap B(\overline{y}%
,\varepsilon)$ such that $w=y-z.$

Take now $\rho\in(0,\min\{(m^{-1}-l)^{-1}\delta,2^{-1}\gamma,2^{-1}%
m\gamma,2^{-1}l^{-1}\gamma,2^{-1}t,4^{-1}r,2^{-1}mt,2^{-1}l^{-1}r\})$ and
define $\beta:=(m^{-1}-l)\rho.$ Also, denote by $\alpha:=\min\{2^{-1}%
\gamma,4^{-1}r\},$ and take $(x,p)\in B(\overline{x},\alpha)\times
B(\overline{p},\alpha).$

If $[\Phi(x,p)-\Psi(x)]\cap B(0,\beta)=\emptyset$ or $0\in\lbrack
\Phi(x,p)-\Psi(x)]\cap B(0,\beta),$ then (\ref{diffix}) trivially holds.
Suppose next that $0\not \in \lbrack\Phi(x,p)-\Psi(x)]\cap B(0,\beta).$ Then,
for every $\xi>0,$ one can find $w_{\xi}\in\lbrack\Phi(x,p)-\Psi(x)]\cap
B(0,\gamma)$ such that%
\begin{equation}
0<\left\Vert w_{\xi}\right\Vert <d(0,[\Phi(x,p)-\Psi(x)]\cap B(0,\beta))+\xi.
\label{ineg}%
\end{equation}

Obviously, $d(0,[\Phi(x,p)-\Psi(x)]\cap B(0,\beta))<(m^{-1}-l)\rho,$ hence for
$\xi>0$ sufficiently small, $d(0,[\Phi(x,p)-\Psi(x)]\cap B(0,\beta
))+\xi<(m^{-1}-l)\rho.$ As consequence, we get from (\ref{ineg}) that%
\begin{equation}
0\in B(w_{\xi},d(0,[\Phi(x,p)-\Psi(x)]\cap B(0,\beta))+\xi)\subset B(w_{\xi
},(m^{-1}-l)\rho)=B(w_{\xi},\beta). \label{0in}%
\end{equation}

\noindent Applying the local sum-stability, one can find $y_{\xi}\in
\Phi(x,p)\cap B({\overline{y},2}^{-1}r)$ and $z_{\xi}\in\Psi(x)\cap
B({\overline{y},2}^{-1}t)$ such that $w_{\xi}=y_{\xi}-z_{\xi}.$ Observe also
that $B(y_{\xi},2^{-1}r)\subset B(\overline{y},r)$ and $B(z_{\xi}%
,2^{-1}t)\subset B(\overline{y},t).$

Denote now $\alpha_{1}:=\beta_{1}:=\alpha_{2}:=\beta_{2}:=2^{-1}\gamma,$
$r_{1}:=2^{-1}t,$ $s_{1}:=2^{-1}t,$ $r_{2}:=4^{-1}r,$ $s_{2}:=2^{-1}r.$
Summarizing,
\begin{align*}
\operatorname*{Gr}\Psi\cap\lbrack D(x,\alpha_{1})\times D(z_{\xi},\beta_{1})]
&  \subset\operatorname*{Gr}\Psi\cap\lbrack D(\overline{x},\gamma)\times
D(\overline{y},\gamma)]\text{ is closed,}\\
\operatorname*{Gr}\Phi_{p}^{-1}\cap\lbrack D(y_{\xi},\beta_{2})\times
D(x,\alpha_{2})]  &  \subset\operatorname*{Gr}\Phi_{p}^{-1}\cap\lbrack
D(\overline{y},\gamma)\times D(\overline{x},\gamma)]\text{ is closed,}\\
\Psi\text{ is }m^{-1}-\text{open at every }(u^{\prime},v^{\prime})  &
\in\operatorname*{Gr}\Psi\cap\lbrack B(x,r_{1})\times B(z_{\xi},s_{1}%
)]\subset\operatorname*{Gr}\Psi\cap\lbrack B(\overline{x},t)\times
B(\overline{y},t)],\\
\Phi_{p}^{-1}\text{ is }l^{-1}-\text{open at every }(y^{\prime},x^{\prime})
&  \in\operatorname*{Gr}\Phi_{p}^{-1}\cap\lbrack B(y_{\xi},s_{2})\times
B(x,r_{2})]\subset\operatorname*{Gr}\Phi_{p}^{-1}\cap\lbrack B(\overline
{y},r)\times B(\overline{x},2^{-1}r)],\\
l^{-1}m^{-1}  &  >1.
\end{align*}
We can apply now Theorem \ref{main_const} for $\Psi,$ $\Phi_{p}^{-1},$
$(x,z_{\xi})\in\operatorname*{Gr}\Psi,$ $(y_{\xi},x)\in\operatorname*{Gr}%
\Phi_{p}^{-1}$ and
\[
\rho_{0}:=(m^{-1}-l)^{-1}(d(0,[\Phi(x,p)-\Psi(x)]\cap B(0,\beta))+\xi
)<\rho<\min\{\alpha_{1},\alpha_{2},m\beta_{1},l^{-1}\beta_{2},r_{1}%
,r_{2},ms_{1},l^{-1}s_{2}\}
\]
to obtain that%
\[
0\in B(z_{\xi}-y_{\xi},d(0,[\Phi(x,p)-\Psi(x)]\cap B(0,\beta))+\xi
)\subset(\Psi-\Phi_{p})(B(x,\rho_{0})).
\]

Using (\ref{0in}), we obtain that $0\in(\Psi-\Phi_{p})(B(x,\rho_{0})),$ so
there exists $\widetilde{x}\in B(x,\rho_{0})$ such that $0\in\Psi
(\widetilde{x})-\Phi(\widetilde{x},p)$ or, equivalently, $\widetilde{x}\in
S(p).$ Hence%
\[
d(x,S(p))\leq d(x,\widetilde{x})<\rho_{0}=(m^{-1}-l)^{-1}(d(0,[\Phi
(x,p)-\Psi(x)]\cap B(0,\beta))+\xi).
\]

Making $\xi\rightarrow0,$ one gets (\ref{diffix}).$\hfill\square$

\bigskip

Notice that in Theorem \ref{fixp} the conclusion reads as a Robinson
regularity of $S$ (see \cite{Rob1976}, \cite{Rob1980}). Remark also that in
order to get (\ref{diffix}) the property of local sum-stability naturally
arises, in contrast to other inequalities involving the Robinson regularity of
$S$ (compare to \cite[Theorem 7]{DonFra2010}, \cite[Theorem 4.4]{DurStr5}).

\bigskip

Finally, Theorem \ref{fixp} gives us the possibility to formulate a result
concerning the calmness of the implicit multifunction $S$ associated to a
parametric variational system.

\begin{thm}
\label{clm_sol}Let $X,Y,P$ be Banach spaces, $F:X\times P\rightrightarrows Y,$
$G:X\rightrightarrows Y$ be two set-valued maps and $(\overline{x}%
,\overline{p},\overline{y})\in X\times P\times Y$ such that $\overline{y}\in
F(\overline{x},\overline{p})$ and $-\overline{y}\in G(\overline{x})$. Suppose
that the following assumptions are satisfied:

(i) $(F,G)$ is locally sum-stable with respect to $x$ uniformly in $p$ around
$((\overline{x},\overline{p}),{\overline{y},}-{\overline{y}});$

(ii) $\operatorname*{Gr}F_{p}$ is locally closed around $(\overline
{x},\overline{y})$ uniformly for $p$ in a neighborhood of $\overline{p};$

(iii) $\operatorname{Gr}G$ is locally closed around $(\overline{x}%
,-\overline{y})$;

(iv) $F$ has the Aubin property with respect to $x,$ uniformly in $p,$ around
$((\overline{x},\overline{p}),\overline{y});$

(v) $F$ is calm with respect to $p,$ uniformly in $x,$ at $((\overline
{x},\overline{p}),\overline{y});$

(vi) $G$ is metrically regular around $(\overline{x},-\overline{y});$

(vii) $\widehat{\operatorname*{lip}}_{x}F((\overline{x},\overline
{p}),\overline{y})\cdot\operatorname*{reg}G(\overline{x},-{\overline{y}})<1.$

Then $S$ is calm at $(\overline{p},\overline{x})$. Moreover, the next relation
is satisfied%
\begin{equation}
\operatorname*{clm}S(\overline{p},\overline{x})\leq\frac{\operatorname*{reg}%
G(\overline{x},-{\overline{y}})\cdot\widehat{\operatorname*{clm}}%
_{p}F((\overline{x},\overline{p}),\overline{y})}{1-\widehat
{\operatorname*{lip}}_{x}F((\overline{x},\overline{p}),\overline{y}%
)\cdot\operatorname*{reg}G(\overline{x},-{\overline{y}})}. \label{clm_S}%
\end{equation}

\end{thm}

\noindent\textbf{Proof. }Take $m>\operatorname*{reg}G(\overline{x}%
,-{\overline{y}})$ and $l>\widehat{\operatorname*{lip}}_{x}F((\overline
{x},\overline{p}),\overline{y})$ such that $m\cdot l<1.$ Apply now Theorem
\ref{fixp} for $F$ and $G$ instead of $\Phi$ and $-\Psi,$ respectively, to get
that there exist $\alpha,\beta>0$ such that for any $(x,p)\in B(\overline
{x},\alpha)\times B(\overline{p},\alpha),$%
\begin{equation}
d(x,S(p))\leq(m^{-1}-l)^{-1}d(0,[F(x,p)+G(x)]\cap B(0,\beta)). \label{mrg}%
\end{equation}

Next, use (v) to get that there exist $\gamma,k>0$ such that for every
$(x,p)\in B({\overline{x},}\gamma)\times B(\overline{p},\gamma),$%
\begin{equation}
e(F_{x}(p)\cap B(\overline{y},\gamma),F_{x}(\overline{p}))\leq kd(p,\overline
{p}). \label{clmFx}%
\end{equation}

From the local sum-stability of $(F,G)$ for $\gamma$ instead of $\varepsilon,$
there is $\delta\in(0,\gamma)$ such that the assertion from Definition
\ref{sum-sta_param} is true. Take now arbitrary $x\in B({\overline{x},}%
\delta),$ $p\in B(\overline{p},\delta)$ and $w\in H(x,p)\cap B(0,\delta),$
where $H$ is given by relation (\ref{H}). Then there exist $y\in F_{x}(p)\cap
B(\overline{y},\gamma)$ and $z\in G(x)\cap B(-\overline{y},\gamma)$ such that
$w=y+z.$ By (\ref{clmFx}), we obtain%
\begin{align*}
d(w,H(x,\overline{p}))  &  =d(y+z,H(x,\overline{p}))\leq d(y+z,F(x,\overline
{p})+z)\\
&  =d(y,F(x,\overline{p}))\leq kd(p,\overline{p}).
\end{align*}

\noindent As $w$ was arbitrarily taken from $H(x,p)\cap B(0,\delta),$ it
follows that $H$ is calm with respect to $p$ uniformly in $x$ at
$((\overline{x},\overline{p}),{0})$. As $k$ can be chosen arbitrarily close to
$\widehat{\operatorname*{clm}}_{p}F((\overline{x},\overline{p}),\overline
{y}),$ we deduce that $\widehat{\operatorname*{clm}}_{p}F((\overline
{x},\overline{p}),\overline{y})\geq\widehat{\operatorname*{clm}}%
_{p}H((\overline{x},\overline{p}),0).$

In conclusion, the relation of the type (\ref{clmS}) follows from (\ref{mrg}),
and $H$ is calm with respect to $p$ uniformly in $x$ at $((\overline
{x},\overline{p}),0),$ so by Theorem \ref{main} we have the conclusion.$\hfill
\square$

\section{Applications to optimization}

We intend to use at-point regularity in solid set-valued optimization
problems. We mention that the incompatibility between efficiency and
around-point regularity generates necessary optimality conditions, as done in
\cite{DurStr1}. Once again, we are interested here to employ the weaker
at-point regularity in the study of multicriteria optimization.

In this section, $X,Y,Z$ are Banach spaces and $K,Q$ are closed convex pointed
cones in $Y$ and $Z$, respectively. As usual, the cone $K$ is proper and gives
a reflexive preorder structure on $Y$ by the equivalence $y_{1}\leq_{K}y_{2}$
if and only if $y_{2}-y_{1}\in K.$ Here the fact that $K$ is proper means that
$K\neq\{0\}$ and $K\neq Y.$

Suppose, in addition, that $K$ is solid, i.e. its topological interior is not
empty ($\operatorname*{int}K\neq\emptyset$).

Then the notion of weak efficiency with respect to the order given by $K$
which we work with is the following.

\begin{df}
\label{defWM}Let $A$ $\subset Y$ be a nonempty subset of $Y.$ A point
$\overline{y}\in A$ is said to be a weak Pareto minimum point of $A$ with
respect to $K$ (we write $\overline{y}\in\operatorname*{WMin}(A,K)$) if
\[
(A-{\overline{y}})\cap(-\operatorname*{int}K)=\emptyset.
\]

\end{df}

A well known feature of solid optimization is that one can use a scalarization
procedure in order to look at a Pareto minimum as a minimum of a scalar
problem (see, for instance, \cite{DurTam}). More precisely, the mechanism is
described by the next result, where $\partial$ denotes the Fenchel
subdifferential of a convex function and $\operatorname{bd}(K)$ denotes the
topological boundary of $K$.

\begin{thm}
\label{scal} For every $e\in\operatorname{int}K,$ the functional
$s_{e}:Y\rightarrow\mathbb{R}$ given by
\begin{equation}
s_{e}(z)=\inf\{\lambda\in\mathbb{R}\mid\lambda e\in z+K\} \label{eq. 1}%
\end{equation}
is continuous, sublinear, strictly-$\operatorname{int}K$-monotone and:

(i) $\partial s_{e}(0)=\{v^{\ast}\in K^{\ast}\mid v^{\ast}(e)=1\}$;

(ii) for every $u\in Y$, $\partial s_{e}(u)\neq\emptyset$ and
\begin{equation}
\partial s_{e}(u)=\{v^{\ast}\in K^{\ast}\mid v^{\ast}(e)=1,v^{\ast}%
(u)=s_{e}(u)\}. \label{subscal}%
\end{equation}

Moreover, $s_{e}$ is $d(e,\operatorname{bd}(K))^{-1}$-Lipschitz. If $A\subset
Y$ is a nonempty set s.t. $A\cap(-\operatorname*{int}K)=\emptyset,$ then
$s_{e}(a)\geq0$ for every $a\in A$.
\end{thm}

For $e\in\operatorname{int}K$ we shall denote $d(e,\operatorname{bd}(K))^{-1}$
by $L_{e}$ (the Lipschitz constant for $s_{e}$).

Now, consider a single-valued map $f:X\rightarrow Y$, a set-valued map
$G:X\rightrightarrows Z$ and the vectorial problem
\[
(P)\hspace{0.4cm}\text{minimize }f(x),\text{ subject to }x\in X,\text{ }0\in
G(x)+Q.
\]

Naturally, a point $\overline{x}\in X$ is called a weak minimum point for
$(P)$ if $f(\overline{x})\in\operatorname{WMin}(f(G^{-1}(-Q)),K)$.

It is well-known that the epigraphical set-valued maps associated to both $f$
and $G$ play an important role in the study of vectorial problem $(P):$ see,
for instance, \cite{BaoMor}, \cite{DurStr1}.

Consider the special type of epigraphical set-valued map associated to $G$ as
$\mathcal{E}_{G}:X\times Z\rightrightarrows Z$ given by%
\[
\mathcal{E}_{G}(x,q):=\left\{
\begin{array}
[c]{ll}%
G(x)+q, & \text{if $q\in Q,$}\\
\emptyset, & \text{otherwise.}%
\end{array}
\right.
\]
This epigraphical multifunction was successfully used in \cite{DNS2011} in
order to give necessary optimality conditions for set-valued optimization
problems without constraints. At this point, just observe that if $G$ is
closed-graph, then the associated multifunction $\mathcal{E}_{G}$ is a
closed-graph multifunction too.

As always when one looks after necessary optimality conditions, one tries to
get some (generalized) Lagrange multipliers in the normal form, i.e. to assure
that the multiplier corresponding to the objective map is non-zero. Of course,
this requires a qualification condition on the constraint system. The problem
$(P)$ was treated in \cite{DurStr3} under metric regularity assumptions. Now
we want to weaken the constraint qualification conditions and for this one
needs the following result which employs a penalization technique for our
problem similar to that in \cite{YY}.

\begin{thm}
\label{thm_scal}Suppose that $\overline{x}\in G^{-1}(-Q)$ is a weak Pareto
minimum point for $(P)$. Fix and denote by $\overline{q}$ a point in $Q$ with
$(\overline{x},\overline{q},0)\in\operatorname*{Gr}\mathcal{E}_{G}.$ If $f$ is
$L-$Lipschitz $(L>0)$ and $\mathcal{E}_{G}$ is metrically subregular at
$((\overline{x},\overline{q}),0)$ (with a constant smaller than $M>0$) then,
for every $e\in\operatorname*{int}K,$ $(\overline{x},\overline{q},0)$ is a
local minimum point for the scalar function%
\[
(x,q,z)\longmapsto s_{e}\circ(f(x)-f(\overline{x}))+LL_{e}M\left\Vert
z\right\Vert
\]
under the constraint $(x,q,z)\in\operatorname*{Gr}\mathcal{E}_{G}.$
\end{thm}

\noindent\textbf{Proof.} Since $f(\overline{x})\in\operatorname{WMin}%
(f(G^{-1}(-Q)),K)$ one has that%

\[
\left[  f(G^{-1}(-Q))-f(\overline{x})\right]  \cap(-\operatorname*{int}%
K)=\emptyset.
\]
Then, following the last assertion in Theorem \ref{scal}, for every
$e\in\operatorname*{int}K$ and $x\in G^{-1}(-Q)$
\begin{equation}
s_{e}(f(x)-f(\overline{x}))\geq0. \label{scal_poz}%
\end{equation}
Now, keeping in mind that the metrical subregularity of $\mathcal{E}_{G}$ at
$((\overline{x},\overline{q}),0)$ is equivalent to the calmness of
$\mathcal{E}_{G}^{-1}$ at $(0,(\overline{x},\overline{q})),$ there exists
$U\in\mathcal{V}(\overline{x}),$ $W\in\mathcal{V}(\overline{q})$ and
$V\in\mathcal{V}(0)$ such that for every $v\in V$
\begin{equation}
\mathcal{E}_{G}^{-1}(v)\cap(U\times W)\subset\mathcal{E}_{G}^{-1}%
(0)+M\left\Vert v\right\Vert D_{X\times Z}(0,1). \label{calm_epi}%
\end{equation}
Let $(x,q,z)\in\operatorname*{Gr}\mathcal{E}_{G}\cap(U\times W\times V).$ From
(\ref{calm_epi}), there is $(x_{0},q_{0})\in\mathcal{E}_{G}^{-1}(0)$ with
\[
\left\Vert (x,q)-(x_{0},q_{0})\right\Vert \leq M\left\Vert z\right\Vert .
\]
Note that, in particular, $0\in G(x_{0})+q\subset G(x_{0})+Q,$ hence $x_{0}$
is a feasible point, and
\[
\left\Vert x-x_{0}\right\Vert \leq M\left\Vert z\right\Vert .
\]
Now one can write, taking firstly into account (\ref{scal_poz}),
\begin{align*}
0  &  =s_{e}\circ(f(\overline{x})-f(\overline{x}))+LL_{e}M\left\Vert
0\right\Vert \\
&  \leq s_{e}(f(x_{0})-f(\overline{x}))=s_{e}(f(x)-f(\overline{x}%
))+s_{e}(f(x_{0})-f(\overline{x}))-s_{e}(f(x)-f(\overline{x}))\\
&  \leq s_{e}(f(x)-f(\overline{x}))+L_{e}\left\Vert f(x)-f(x_{0})\right\Vert
\leq s_{e}(f(x)-f(\overline{x}))+L_{e}L\left\Vert x-x_{0}\right\Vert \\
&  \leq s_{e}(f(x)-f(\overline{x}))+L_{e}LM\left\Vert z\right\Vert .
\end{align*}
This proves the assertion in conclusion.\hfill$\square$

\medskip

In order to give dual conditions for metric subregularity of $\mathcal{E}_{G}$
we need to recall the mechanism of error bounds of a system given by
inequality constraints.

Let $X$ be a normed vector space and let $f:X\rightarrow\mathbb{R}%
\cup\{+\infty\}$ be a given function. We set
\begin{equation}
S:=\{x\in X:\quad f(x)\leq0\}. \label{1}%
\end{equation}
One denotes the quantity $\max\{f(x),0\}$ by $[f(x)]_{+}.$ We say that the
system (\ref{1}) admits an error bound if there exists a real $c>0$ such that
\begin{equation}
d(x,S)\leq c[f(x)]_{+}\quad\text{for all }x\in X. \label{2}%
\end{equation}
For $x_{0}\in S$, we say that the system (\ref{1}) has an error bound at
$x_{0}$, if there exists a real $c>0$ such that relation (\ref{2}) is
satisfied for all $x$ around $x_{0}$.

Here and in what follows the convention $0\cdot(+\infty)=0$ is used.

The following result will be useful in the sequel. Notice that the parametric
case is studied in \cite{NTT}.

\begin{thm}
\label{Err_bd}(\cite[Corollary 2.3]{SiamNT}) Let $X$ be a Banach space, and
let $f:X\rightarrow\mathbb{R}\cup\{+\infty\}$ be a lower semicontinuous
function. Let $\overline{x}\in S,$ $\tau\in(0,+\infty)$ and $\eta\in
(0,+\infty)$ be given. Consider the following statements:

(i) $d(x,S)\leq\tau\lbrack f(x)]_{+},$\quad for all $x\in B(\overline{x}%
,\eta/2).$

(ii) For each $x\in B(\overline{x},\eta)\backslash S$ and for any
$\varepsilon>0$, there exists $z\in X$ such that%
\begin{equation}
0<d(x,z)<(\tau+\varepsilon)(f(x)-[f(z)]_{+}). \label{ineq1}%
\end{equation}

(iii) For each $x\in B(\overline{x},\eta)\setminus S$ and for any
$\varepsilon>0$, there exists $z\in X$ with $f(z)\geq0$ such that
(\ref{ineq1}) holds.

(iv) For each $x\in B(\overline{x},\eta)\setminus S$ and for any
$\varepsilon>0$, there exists $z\in X$ with $f(z)>0$ such that (\ref{ineq1}) holds.

Then, one has $(iv)\Rightarrow(iii)\Rightarrow(ii)\Rightarrow(i)$. Conversely,
if $(i)$ holds, then $(ii)$ holds with $\eta/2$ instead of $\eta$. In
addition, if $f$ is a continuous function, then the three statements $(ii)$,
$(iii)$, and $(iv)$ are equivalent.
\end{thm}

\bigskip

The strong slope $|\nabla f|(x)$ of a lower semicontinuous function $f$ at
$x\in\operatorname*{dom}f:=\{u\in X:f(u)<+\infty\}$ is the quantity defined by
$|\nabla f|(x)=0$ if $x$ is a local minimum of $f,$ and by
\[
|\nabla f|(x)=\limsup_{y\rightarrow x,y\neq x}\frac{f(x)-f(y)}{d(x,y)},
\]
otherwise. For $x\notin\operatorname*{dom}f,$ we set $|\nabla f|(x)=+\infty$
(see \cite{DMT}, \cite{AC}). From Theorem \ref{Err_bd} we get the following
result. In order to clarify the ideas, we present its proof.

\begin{cor}
\label{Slope1}Let $X$ be a Banach space and let $f:X\rightarrow\mathbb{R}%
\cup\{+\infty\}$ be a lower semicontinuous function. Let $\gamma\in
(0,+\infty)$ and $\overline{x}\in S$ be given. If there exist a neighborhood
$V$ of $\overline{x}$ and a real $m>0$ such that $|\nabla f|(x)\geq m$ for all
$x\in V$ with $f(x)\in(0,\gamma)$ then there exists a neighborhood $V_{1}$of
$\overline{x}$ such that
\begin{equation}
md(x,S)\leq\lbrack f(x)]_{+}\text{\quad for all }x\in V_{1}. \label{errbd}%
\end{equation}

\end{cor}

\noindent\textbf{Proof. }Without loosing the generality, suppose that
$V=B(\overline{x},\eta)$ with $\eta>0.$ Observe first that if $x\in V$ is such
that $f(x)\leq0,$ then (\ref{errbd}) trivially holds.

Fix now $x\in V$ such that $f(x)\in(0,\gamma).$ But this imply, on one hand,
that $|\nabla f|(x)<+\infty,$ and on the other hand, that there is $\alpha>0$
such that $f(x)>\alpha>0.$ Using now the lower semicontinuity of $f,$ we can
find a neighborhood $U$ of $x$ such that, for every $u\in U,$ $f(u)>\alpha>0.$

Using now the assumptions of the corollary, we also know that $|\nabla
f|(x)\geq m>0,$ so $x$ cannot be a local minimum of the function $f.$ In
conclusion, by the definition of the strong slope, for every $\varepsilon>0,$
and for every $N\in\mathcal{V}(x),$ there exists $z\in N,$ $z\not =x$ such
that%
\begin{equation}
(m-\varepsilon)^{-1}(f(x)-f(z))\geq d(z,x)>0. \label{est}%
\end{equation}

By taking $N$ sufficiently small, such that $N\subset U,$ we obtain that
$0<f(z)=[f(z)]_{+}.$ In conclusion, for every $x\in B(\overline{x}%
,\eta)\setminus S$ and every $\varepsilon>0,$ we have found $z\in X$ with
$f(z)>0$ such that%
\[
(m-\varepsilon)^{-1}(f(x)-[f(z)]_{+})\geq d(z,x)>0.
\]

The conclusion now follows from the implication $(iv)\Rightarrow(i)$ of the
previous theorem.\hfill$\square$

\medskip

Coming back to our epigraphical set-valued map $\mathcal{E}_{G},$ define the
lower semicontinuous application $\varphi_{\mathcal{E}_{G}}:X\times Z\times
Z\rightarrow\mathbb{R\cup\{+\infty\}}$ by%

\[
\varphi_{\mathcal{E}_{G}}((x,q),z):=\liminf_{(u,r,v)\rightarrow(x,q,z)}%
d(v,\mathcal{E}_{G}(u,r))
\]
which, as shown in \cite{DNS2011}, takes the form%
\[
\varphi_{\mathcal{E}_{G}}((x,q),z)=\left\{
\begin{array}
[c]{ll}%
\liminf\limits_{u\rightarrow x}d(z,G(u)+q), & \text{if }q\text{$\in Q,$}\\
+\infty, & \text{otherwise.}%
\end{array}
\right.
\]

\bigskip

Note that, if $G$ is a closed-graph multifunction, then for every $z\in Z$,
\[
\{(x,q)\in X\times Z:y\in G(x)+q,q\in Q\}=\{(x,q)\in X\times Z:\varphi
_{\mathcal{E}_{G}}((x,q),z)=0\}.
\]

\begin{thm}
\label{mreg_EF}Let $X$ be a Banach space, let $Z$ be a normed vector space and
let $G:X\rightrightarrows Z$ be a closed-graph multifunction. Suppose that
$(\overline{x},\overline{q},\overline{z})\in X\times Z\times Z$ is such that
$\overline{z}\in G(\overline{x})+\overline{q}$ and $\overline{q}\in Q.$ Let
$m>0$ be given. If there exist a neighborhood $U\times V$ of $(\overline
{x},\overline{q})$ and a real $\gamma>0$ such that
\begin{equation}
|\nabla\varphi_{\mathcal{E}_{G}}((\cdot,\cdot),\overline{z})|(x,q)\geq
m\quad\text{for all}\;(x,q)\in U\times V\ \text{with}\;\varphi_{\mathcal{E}%
_{G}}((x,q),\overline{z})\in(0,\gamma), \label{fiG}%
\end{equation}
then there exists a neighborhood $\tilde{U}\times\tilde{V}\ $of $(\overline
{x},\overline{q})$ such that
\begin{equation}
d((x,q),\mathcal{E}_{G}^{-1}(\overline{z}))\leq\frac{1}{m}d(\overline
{z},G(x)+q)\quad\forall(x,q)\in\tilde{U}\times\lbrack\tilde{V}\cap Q].
\label{mregG}%
\end{equation}

In other words, $\mathcal{E}_{G}$ is metrically subregular at $((\overline
{x},\overline{q}),\overline{z}),$ hence $\mathcal{E}_{G}$ is linearly
pseudo-open at $((\overline{x},\overline{q}),\overline{z}).$
\end{thm}

\noindent\textbf{Proof. }We apply Corollary \ref{Slope1} for the function%
\[
(x,q)\mapsto\varphi_{\mathcal{E}_{G}}((x,q),\overline{z}).
\]
Note that in this case
\begin{align*}
S  &  =\{(x,q)\in X\times Z:\varphi_{\mathcal{E}_{G}}((x,q),\overline
{z})=0\}\\
&  =\{(x,q)\in X\times Z:\overline{z}\in G(x)+q,\text{ }q\in Q\}\\
&  =\mathcal{E}_{G}^{-1}(\overline{z})
\end{align*}
and obviously $(\overline{x},\overline{q})\in S.$ Since the assumptions of
Corollary \ref{Slope1} are fulfilled, we can find a neighborhood $\tilde
{U}\times\tilde{V}\ $of $(\overline{x},\overline{q})$ such that%
\[
md((x,q),S)\leq\lbrack\varphi_{\mathcal{E}_{G}}((x,q),\overline{z}%
)]_{+}=\varphi_{\mathcal{E}_{G}}((x,q),\overline{z}),\text{\quad for all
}(x,q)\in\tilde{U}\times\lbrack\tilde{V}\cap Q].
\]
Since $\varphi_{\mathcal{E}_{G}}((x,q),\overline{z})\leq d(\overline
{z},G(x)+q)$ the conclusion follows.$\hfill\square$

\medskip

The main tools for our study are the Mordukhovich's generalized
differentiation objects which are the basis of many effective techniques in
variational analysis. We recall the most important facts we need in this paper.

\bigskip

Let $X$ be a normed vector space, $S$ be a non-empty subset of $X$ and let
$x\in S.$ The Fr\'{e}chet normal cone to $S$ at $x$ is%
\begin{equation}
\widehat{N}(S,x):=\left\{  x^{\ast}\in X^{\ast}\mid\underset{u\overset
{S}{\rightarrow}x}{\lim\sup}\frac{x^{\ast}(u-x)}{\left\Vert u-x\right\Vert
}\leq0\right\}  . \label{eps-no}%
\end{equation}

If $X$ is an Asplund space (i.e. a Banach space where every convex continuous
function is generically Fr\'{e}chet differentiable), the basic (or limiting,
or Mordukhovich) normal cone to $S$ at $\overline{x}$ is:%
\[
N(S,\overline{x})=\{x^{\ast}\in X^{\ast}\mid\exists x_{n}\overset
{S}{\rightarrow}\overline{x},x_{n}^{\ast}\overset{w^{\ast}}{\rightarrow
}x^{\ast},x_{n}^{\ast}\in\widehat{N}(S,x_{n}),\forall n\in\mathbb{N}\}.
\]

Let $f:X\rightarrow\overline{\mathbb{R}},$ finite at $\overline{x}\in X;$ the
Fr\'{e}chet subdifferential of $f$ at $\overline{x}$ is the set
\[
\widehat{\partial}f(\overline{x}):=\{x^{\ast}\in X^{\ast}\mid(x^{\ast}%
,-1)\in\widehat{N}(\operatorname*{epi}f,(\overline{x},f(\overline{x})))\}
\]
and, if $X$ is Asplund, the basic (or limiting, or Mordukhovich)
subdifferential of $f$ at $\overline{x}$ is%

\[
\partial f(\overline{x}):=\{x^{\ast}\in X^{\ast}\mid(x^{\ast},-1)\in
N(\operatorname*{epi}f,(\overline{x},f(\overline{x})))\},
\]
where $\operatorname*{epi}f$ denotes the epigraph of $f.$ On Asplund spaces
one has
\[
\partial f(\overline{x})=\limsup_{x\overset{f}{\rightarrow}\overline{x}%
}\widehat{\partial}f(x),
\]
and, in particular, $\widehat{\partial}f(\overline{x})\subset\partial
f(\overline{x}).$ Note that a generalized Fermat rule holds: if $\overline{x}$
is a local minimum point for $f$ then $0\in\widehat{\partial}f(\overline{x}).$
If $f$ is convex, then both these subdifferentials do coincide with the
classical Fenchel subdifferential. If $\delta_{\Omega}$ denotes the indicator
function associated with a nonempty set $\Omega\subset X$ (i.e. $\delta
_{\Omega}(x)=0$ if $x\in\Omega,$ $\delta_{\Omega}(x)=\infty$ if $x\notin
\Omega$), then for any $\overline{x}\in\Omega,$ $\widehat{\partial}%
\delta_{\Omega}(\overline{x})=\widehat{N}(\Omega,\overline{x})$ and
$\partial\delta_{\Omega}(\overline{x})=N(\Omega,\overline{x}).$ Let
$\Omega\subset X$ be a nonempty set and take $\overline{x}\in\Omega;$ then one
has:%
\begin{equation}
\widehat{\partial}d(\cdot,\Omega)(\overline{x})=\widehat{N}(\Omega
,\overline{x})\cap D_{X^{\ast}},\text{ }\widehat{N}(\Omega,\overline
{x})=\underset{\lambda>0}{%
{\displaystyle\bigcup}
}\lambda\widehat{\partial}d(\cdot,\Omega)(\overline{x}). \label{s_dist}%
\end{equation}

The basic subdifferential satisfies a robust sum rule (see \cite[Theorem
3.36]{Mor2006}): if $X$ is Asplund, $f_{1}$ is Lipschitz around $\overline{x}$
and $f_{2}$ is lower semicontinuous around this point, then
\begin{equation}
\partial(f_{1}+f_{2})(\overline{x})\subset\partial f_{1}(\overline
{x})+\partial f_{2}(\overline{x}). \label{r_bas_calc}%
\end{equation}

A function $f:X\rightarrow Y$ is said to be strictly Lipschitz at
$\overline{x}$ if it is locally Lipschitzian around this point and there
exists a neighborhood $V$ of the origin in $X$ s.t. the sequence $(t_{k}%
^{-1}(f(x_{k}+t_{k}v)-f(x_{k})))_{k\in\mathbb{N}}$ contains a norm convergent
subsequence whenever $v\in V,x_{k}\rightarrow\overline{x},$ $t_{k}%
\downarrow0.$

Suppose that $X,Y$ are Asplund spaces. Let $f:X\rightarrow Y$ and
$\varphi:Y\rightarrow\mathbb{R}$ s.t. $f$ is strictly Lipschitz at
$\overline{x}\in X$ and $\varphi$ is Lipschitz around $f(\overline{x});$ then
\begin{equation}
\partial(\varphi\circ f)(\overline{x})\subset%
{\displaystyle\bigcup\limits_{y^{\ast}\in\partial\varphi(f(\overline{x}))}}
\partial(y^{\ast}\circ f)(\overline{x}). \label{rule_chain}%
\end{equation}

Let $F:X\rightrightarrows Y$ be a set-valued map and $(\overline{x}%
,\overline{y})\in\operatorname*{Gr}F.$ Then the Fr\'{e}chet coderivative at
$(\overline{x},\overline{y})$ is the set-valued map $\widehat{D}^{\ast
}F(\overline{x},\overline{y}):Y^{\ast}\rightrightarrows X^{\ast}$ given by
\[
\widehat{D}^{\ast}F(\overline{x},\overline{y})(y^{\ast}):=\{x^{\ast}\in
X^{\ast}\mid(x^{\ast},-y^{\ast})\in\widehat{N}(\operatorname{Gr}%
F,(\overline{x},\overline{y}))\}.
\]

Similarly, on Asplund spaces, the normal coderivative of $F$ at $(\overline
{x},\overline{y})$ is the set-valued map $D^{\ast}F(\overline{x},\overline
{y}):Y^{\ast}\rightrightarrows X^{\ast}$ given by
\[
D^{\ast}F(\overline{x},\overline{y})(y^{\ast}):=\{x^{\ast}\in X^{\ast}%
\mid(x^{\ast},-y^{\ast})\in N(\operatorname{Gr}F,(\overline{x},\overline
{y}))\}.
\]

\bigskip

Now, the next theorem is a reformulation of the main result in \cite{DNS2011},
where by $Q^{\ast}$ we denote the positive polar cone of $Q$.

\begin{thm}
\label{EstimSlopeEF} Let $X,Z$ be Banach spaces and $G$ be closed-graph. Then
for each $(x,q_{0},z)\in X\times Q\times Z$ with $z\notin G(x)+q_{0},$ one
has
\begin{equation}
|\nabla\varphi_{\mathcal{E}_{G}}((\cdot,\cdot),z)|(x,q_{0})\geq\lim
_{\rho\downarrow0}\left\{  \inf\left\{  \Vert x^{\ast}\Vert:\;%
\begin{array}
[c]{l}%
(u,v)\in\operatorname*{Gr}G,\;u\in B(x,\rho),x^{\ast}\in\widehat{D}^{\ast
}G(u,v)(z^{\ast}+v^{\ast}),\\
z^{\ast}\in Q^{\ast}\cap S_{Z^{\ast}},\text{ }v^{\ast}\in\rho B_{Z^{\ast}},\\
d(z,G(u)+q_{0})\leq\varphi_{\mathcal{E}_{G}}((x,q_{0}),z)+\rho,\\
\Vert z-v-q_{0}\Vert\leq d(z,G(u)+q_{0})+\rho,\\
|\langle z^{\ast}+v^{\ast},z-v-q_{0}\rangle-d(z,G(u)+q_{0})|<\rho
\end{array}
\right\}  \right\}  . \label{est_ssl}%
\end{equation}

\end{thm}

\medskip

Now, putting together Theorems \ref{mreg_EF} and \ref{EstimSlopeEF} one gets a
metric subregularity sufficient condition in terms of dual objects.

\begin{thm}
\label{openEF}Let $X,Z$ be Banach spaces and $G$ be closed-graph. Suppose that
$(\overline{x},\overline{q},\overline{z})\in X\times Z\times Z$ is such that
$\overline{z}\in G(\overline{x})+\overline{q}$ and $\overline{q}\in Q.$ Let
$m>0$ be given. If there exist a neighborhood $U\times V$ of $(\overline
{x},\overline{q})$ and a real $\gamma>0$ such that for every $(x,q)\in
U\times\lbrack V\cap Q]$ with $\overline{z}\not \in G(x)+q,$%
\begin{equation}
m\leq\lim_{\rho\downarrow0}\left\{  \inf\left\{  \Vert x^{\ast}\Vert:\;%
\begin{array}
[c]{l}%
(u,v)\in\operatorname*{Gr}G,\;u\in B(x,\rho),x^{\ast}\in\widehat{D}^{\ast
}G(u,v)(z^{\ast}+v^{\ast}),\\
z^{\ast}\in Q^{\ast}\cap S_{Z^{\ast}},\text{ }v^{\ast}\in\rho B_{Z^{\ast}},\\
d(\overline{z},G(u)+q)\leq\gamma+\rho,\\
\Vert\overline{z}-v-q\Vert\leq d(\overline{z},G(u)+q)+\rho,\\
|\langle z^{\ast}+v^{\ast},\overline{z}-v-q\rangle-d(\overline{z}%
,G(u)+q)|<\rho
\end{array}
\right\}  \right\}  \label{m}%
\end{equation}
then $\mathcal{E}_{G}$ is metrically subregular at $((\overline{x}%
,\overline{q}),\overline{z}).$
\end{thm}

\noindent\textbf{Proof. }Take arbitrary $(x,q)\in U\times V\ $with$\;\varphi
_{\mathcal{E}_{G}}((x,q),\overline{z})\in(0,\gamma).$ Then $q\in Q,$ because
otherwise $\varphi_{\mathcal{E}_{G}}((x,q),\overline{z})=+\infty.$ Also,
$\overline{z}\not \in G(x)+q,$ because otherwise $\varphi_{\mathcal{E}_{G}%
}((x,q),\overline{z})=0.$ Now, since the set from (\ref{est_ssl}) (where the
infimum is taken) is smaller than the corresponding set from (\ref{m}), we get
that $|\nabla\varphi_{\mathcal{E}_{G}}((\cdot,\cdot),\overline{z})|(x,q)\geq
m.$ But this means, in virtue of Theorem \ref{mreg_EF}, exactly the
conclusion.\hfill$\square$

\medskip

In some cases, the condition (\ref{m}) could be simplified, as the next
corollary shows.

\begin{cor}
\label{th_opt}Let $X,Z$ be Asplund spaces and $G$ be closed-graph. Suppose
that $(\overline{x},\overline{z})\in X\times Z$ is such that $\overline{z}\in
G(\overline{x}).$ If there exist $r,c>0$ such that, for every $(x,q,z)\in
B(\overline{x},r)\times B(0,r)\times B(\overline{z},r),$ $(x,z)\in
\operatorname*{Gr}G,$ $(x,q,\overline{z})\notin\operatorname*{Gr}%
\mathcal{E}_{G},$ $z^{\ast}\in Q^{\ast}\cap S_{Z^{\ast}},v^{\ast}%
\in2cB_{Y^{\ast}},x^{\ast}\in D^{\ast}G(x,z)(z^{\ast}+v^{\ast}),$%
\[
c\left\Vert z^{\ast}+v^{\ast}\right\Vert \leq\left\Vert x^{\ast}\right\Vert ,
\]

\noindent then $\mathcal{E}_{G}$ linearly pseudo-open at $((\overline
{x},0),\overline{z}).$
\end{cor}

\noindent\textbf{Proof. }Fix $m:=2^{-1}c>0,$ $\gamma:=2^{-1}r>0,$
$U:=B(\overline{x},2^{-1}r),$ $V:=B(0,4^{-1}r).$

Suppose first that for arbitrary $(x,q)\in U\times\lbrack V\cap Q],$ one has
that $\overline{z}\in G(x)+q=\mathcal{E}_{G}(x,q).$ It follows in particular
that for every $(x,q)\in B(\overline{x},2^{-1}r)\times\lbrack B(0,r)\cap Q],$
one has that $\overline{z}\in\mathcal{E}_{G}(x,q).$ Consequently, for every
$\tau\in(0,2^{-1}r),$ one has that%
\[
\overline{z}\in\mathcal{E}_{G}(x,q)\subset\mathcal{E}_{G}(B(x,\tau
),B(q,\tau)).
\]

As the previous relation is true for every $(x,q)\in U\times\lbrack V\cap Q]$
and every $\tau\in(0,2^{-1}r),$ we conclude that $\mathcal{E}_{G}$ is linearly
pseudo-open at $((\overline{x},0),\overline{z}).$

Suppose now that there exists $(x,q)\in U\times\lbrack V\cap Q]$ such that
$\overline{z}\not \in G(x)+q$, i.e. $(x,q,\overline{z})\notin
\operatorname*{Gr}\mathcal{E}_{G}.$ Choose $\rho\in(0,\min\{8^{-1}%
r,2c,2^{-1}\}).$ Consider $(u,v)\in\operatorname*{Gr}G,\;u\in B(x,\rho),$
$z^{\ast}\in Q^{\ast}\cap S_{Z^{\ast}},$ $v^{\ast}\in\rho B_{Z^{\ast}},$
$x^{\ast}\in\widehat{D}^{\ast}G(u,v)(z^{\ast}+v^{\ast}),$ $\Vert\overline
{z}-v-q\Vert\leq d(\overline{z},G(u)+q)+\rho,$ $d(\overline{z},G(u)+q)\leq
\gamma+\rho,$ $|\langle z^{\ast}+v^{\ast},\overline{z}-v-q\rangle
-d(\overline{z},G(u)+q)|<\rho.$ Then%
\begin{align*}
\left\Vert u-\overline{x}\right\Vert  &  \leq\left\Vert u-x\right\Vert
+\left\Vert x-\overline{x}\right\Vert <\rho+2^{-1}r<r,\\
\left\Vert v-\overline{z}\right\Vert  &  \leq\left\Vert v+q-\overline
{z}\right\Vert +\left\Vert q\right\Vert <d(\overline{z},G(u)+q)+\rho+4^{-1}r\\
&  \leq\gamma+2\rho+4^{-1}r<2^{-1}r+4^{-1}r+4^{-1}r=r,\\
\left\Vert v^{\ast}\right\Vert  &  <\rho<2c.
\end{align*}

One can use now the hypothesis from the statement of the Corollary to get that%
\[
\left\Vert x^{\ast}\right\Vert \geq c\left\Vert z^{\ast}+v^{\ast}\right\Vert
\geq c(\left\Vert z^{\ast}\right\Vert -\left\Vert v^{\ast}\right\Vert )\geq
c(1-\rho)\geq m.
\]

\noindent Using now Theorem \ref{openEF}, one gets the conclusion.\hfill
$\square$

\medskip

Summing up, Theorem \ref{thm_scal} gives a scalarization method for $(P)$
under at-point assumption of the epigraphical multifunction associated to the
constraints, while Theorem \ref{mreg_EF} provides sufficient conditions for
the fulfilment of this assumption. Finally, taking advantage of the power of
Mordukhovich subdifferential calculus we present necessary optimality
conditions for $(P)$ in terms of generalized differentiation.

\begin{thm}
Take $X,Y,Z$ as Asplund spaces and $\overline{x}\in G^{-1}(-Q)$ as a weak
Pareto minimum point for $(P)$. Fix $\overline{q}\in Q$ such that
$(\overline{x},\overline{q},0)\in\operatorname*{Gr}\mathcal{E}_{G}.$ Suppose that

(i) $f$ is $L-$Lipschitz $(L>0)$ and strictly Lipschitz at $\overline{x};$

(ii) $G$ is closed-graph;

(iii) $\mathcal{E}_{G}$ is metrically subregular at $((\overline{x}%
,\overline{q}),0)$ (with a constant smaller than $M>0$).

Then, for every $e\in\operatorname*{int}K,$ there exist $y^{\ast}\in K^{\ast
},$ $y^{\ast}(e)=1,$ $z^{\ast}\in Z^{\ast},$ $\left\Vert z^{\ast}\right\Vert
\leq LL_{e}M$ s.t.
\[
(0,0)\in\partial(y^{\ast}\circ f)(\overline{x})\times\{0\}+D^{\ast}%
\mathcal{E}_{G}(\overline{x},\overline{q},0)(z^{\ast}).
\]

\end{thm}

\noindent\textbf{Proof.} Consider $e\in\operatorname*{int}K.$ Taking into
account Theorem \ref{thm_scal}, $(\overline{x},\overline{q},0)$ is a local
minimum point for the scalar function%
\[
(x,q,z)\longmapsto s_{e}\circ(f(x)-f(\overline{x}))+LL_{e}M\left\Vert
z\right\Vert
\]
under the constraint $(x,q,z)\in\operatorname*{Gr}\mathcal{E}_{G}.$ As usual,
this means that $(\overline{x},\overline{q},0)$ is a local minimum point for
the unconstrained scalar problem%
\[
\min\left[  s_{e}\circ(f(\cdot)-f(\overline{x}))+LL_{e}M\left\Vert
\cdot\right\Vert +\delta_{\operatorname*{Gr}\mathcal{E}_{G}}(\cdot,\cdot
,\cdot)\right]  .
\]
Consequently,
\[
(0,0,0)\in\partial\left[  s_{e}\circ(f(\cdot)-f(\overline{x}))+LL_{e}%
M\left\Vert \cdot\right\Vert +\delta_{\operatorname*{Gr}\mathcal{E}_{G}}%
(\cdot,\cdot,\cdot)\right]  (\overline{x},\overline{q},0).
\]
Since under our assumption we can apply the exact sum rule (\ref{r_bas_calc})
we get (by means of some obvious calculations):%
\[
(0,0,0)\in\partial s_{e}\circ(f(\cdot)-f(\overline{x}))(\overline{x}%
)\times\{0\}\times\{0\}+LL_{e}M\left[  \{0\}\times\{0\}\times D_{Z^{\ast}%
}(0,1)\right]  +\partial\delta_{\operatorname*{Gr}\mathcal{E}_{G}}%
(\overline{x},\overline{q},0).
\]
The hypothesis $(i)$ gives us the right to use the chain rule
(\ref{rule_chain}): there exists $y^{\ast}\in\partial s_{e}(0)$ s.t.%
\[
(0,0,0)\in\partial(y^{\ast}\circ f)(\overline{x})\times\{0\}\times
\{0\}+LL_{e}M\left[  \{0\}\times\{0\}\times D_{Z^{\ast}}(0,1)\right]
+N(\operatorname*{Gr}\mathcal{E}_{G},(\overline{x},\overline{q},0)).
\]
Then there exist $u^{\ast}\in\partial(y^{\ast}\circ f)(\overline{x}),$
$z^{\ast}\in D_{Z^{\ast}}(0,LL_{e}M)$ with
\[
(-u^{\ast},0,-z^{\ast})\in N(\operatorname*{Gr}\mathcal{E}_{G},(\overline
{x},\overline{q},0)),
\]
i.e.
\[
(-u^{\ast},0)\in D^{\ast}\mathcal{E}_{G}(\overline{x},\overline{q},0)(z^{\ast
})
\]
which achieves the proof.\hfill$\square$

\end{document}